\documentclass[12pt]{amsart}
\usepackage{latexsym, amssymb, amsmath, amscd}
 \usepackage{eucal}

    \newtheorem{rema}{Remark}[section]
    \newtheorem{propo}[rema]{Proposition}

   \newtheorem{theo}[rema]{Theorem}
   \newtheorem{def-theo}[rema]{Definition-Theorem}
 \newtheorem{conj}[rema]{Conjecture}
   \newtheorem{defi}[rema]{Definition}
    \newtheorem{lemma}[rema]{Lemma}
    \newtheorem{corol}[rema]{Corollary}
     
  \newtheorem{rmk}[rema]{Remark}

    \newtheorem{prob}[rema]{Problem}

	\newcommand{\nno}{\nonumber}

	\newcommand{\p}{\partial}

 \newcommand{\pf}{{\it Proof:}\hspace{2ex}}
 
 \newcommand{\epfv}{\hspace{1em}$\Box$\vspace{1em}}

\newcommand{\bC}{{\mathbb C}}
\newcommand{\bZ}{{\mathbb Z}}

\newcommand{\bQ}{{\mathbb Q}}
\newcommand{\bN}{{\mathbb N}}

\newcommand{\cS}{{\mathcal S}}

\newcommand{\cD}{{\mathcal D}}

\newcommand{\cC}{{\mathcal C}}

\newcommand{\cNcs}{{${\mathcal N}$CS} }
\newcommand{\cNsf}{{{\mathcal N}Sym}}

\newcommand{\Oft}{ \Omega_{F_t}  }

\newcommand{\cSft}{ {\mathcal S}_{F_t} }

\newcommand{\cDaz}{{{\mathcal D}^{[\alpha]} \langle  z \rangle}}

\newcommand{\cDzz}{{\mathcal D \langle \langle z \rangle\rangle}}
\newcommand{\cDkzz}{{\mathcal D_K \langle \langle z \rangle\rangle}}
\newcommand{\cDtzz}{{\mathcal D_t \langle \langle z \rangle\rangle}}
\newcommand{\cDrzz}{{\mathcal D er\langle \langle z \rangle\rangle}}
\newcommand{\cDrtzz}{{\mathcal D er_t \langle \langle z \rangle\rangle}}
\newcommand{\cDrkzz}{{\mathcal D er_K \langle \langle z \rangle\rangle}}

\newcommand{\cDrazz}{{\cD er^{[\alpha]}\langle \langle z \rangle \rangle}} 
 
\newcommand{\cDazz}{{\cD^{[\alpha]}\langle \langle z \rangle \rangle}}

\newcommand{\cDrtazz}{{\cD er^{[\alpha]}_t\langle \langle z \rangle \rangle}} 
 
\newcommand{\cDtazz}{{\cD^{[\alpha]}_t\langle \langle z \rangle \rangle}}

\newcommand{\ataz}{{\mathbb A_t^{[\alpha]}\langle \langle z\rangle\rangle}}

\newcommand{\etaz}{{\mathbb E_t^{[\alpha]}\langle \langle z\rangle\rangle}}

\newcommand{\kzz}{{K\langle \langle z \rangle\rangle}}
\newcommand{\ktz}{{K[t]\langle z \rangle}}
\newcommand{\ktzz}{{K[t]\langle \langle z \rangle\rangle}}
\newcommand{\kttzz}{{K[[t]]\langle \langle z \rangle\rangle}}

\newcommand{\BQ}{\begin{eqnarray}}
\newcommand{\EQ}{\end{eqnarray}}
\newcommand{\BQn}{\begin{eqnarray*}}
\newcommand{\EQn}{\end{eqnarray*}}

\newcommand{\lb}{\left[}
\newcommand{\rb}{\right]}
\newcommand{\lp}{\left(}
\newcommand{\rp}{\right)}

\newcommand{\fr}{\frac}
\newcommand{\pz}{\frac{\p}{\p z}}

\title[NCSF's and the Inversion Problem]
{Noncommutative Symmetric Functions and the Inversion Problem}

    \author{Wenhua Zhao}      
    \begin{document}

\begin{abstract}
Let $K$ be any unital commutative $\bQ$-algebra and 
$z=(z_1, z_2, \cdots, z_n)$ commutative 
or noncommutative variables. Let $t$ be a 
formal central parameter and $\kttzz$ 
the formal power series algebra of 
$z$ over $K[[t]]$. In \cite{GTS-II}, 
for each automorphism $F_t(z)=z-H_t(z)$ 
of $\kttzz$ with $H_{t=0}(z)=0$ 
and $o(H(z))\geq 1$, a \cNcs (noncommutative symmetric) system 
(\cite{GTS-I}) $\Oft$ has been constructed. 
Consequently, we get a Hopf algebra homomorphism  
$\cSft: \cNsf \to \cDzz$ from the 
Hopf algebra $\cNsf$ (\cite{G-T}) 
of NCSF's (noncommutative symmetric functions). 
In this paper, 
we first give a list for the identities 
between any two sequences of differential operators 
in the \cNcs  system $\Oft$ by using some 
identities of NCSF's derived in \cite{G-T} 
and the homomorphism $\cSft$. 
Secondly, we apply these identities 
to derive some formulas in terms of differential operator in 
the system $\Oft$ for the Taylor series expansions of 
$u(F_t)$ and $u(F_t^{-1})$ $(u(z)\in \kttzz)$; 
the D-Log and the formal flow of $F_t$ and 
inversion formulas for the inverse map of $F_t$.
Finally, we discuss a connection 
of the well-known Jacobian conjecture 
with NCSF's.
\end{abstract}

   

 \bibliographystyle{alpha}
    \maketitle


\renewcommand{\theequation}{\thesection.\arabic{equation}}
\renewcommand{\therema}{\thesection.\arabic{rema}}
\setcounter{equation}{0}
\setcounter{rema}{0}
\setcounter{section}{0}

\section{\bf Introduction}\label{S1}

Let $K$ be any unital commutative $\bQ$-algebra and 
$z=(z_1, z_2, \cdots, z_n)$ commutative 
or noncommutative variables. Let $t$ be a 
formal central parameter, 
i.e. a formal variable which commutes 
with $z$ and elements of $K$.
To keep notation simple, 
we use the notations 
for noncommutative variables 
uniformly for both commutative 
and noncommutative variables $z$. 
Let $\kzz$ (resp.\,\,$\kttzz$) 
the algebra of formal power series 
in $z$ over $K$ (resp.\,\,$K[[t]]$). 
For any $\alpha \geq 1$, 
let $\cDazz$ be the unital algebra 
generated by 
the differential operators of $\kzz$ 
which increase the degree in 
$z$ by at least $\alpha-1$ 
and $ \ataz $ the group of automorphisms 
$F_t(z)=z-H_t(z)$ 
of $\kttzz$ with $o(H_t(z))\geq \alpha$ 
and $H_{t=0}(z)=0$.
In \cite{GTS-II}, for each automorphism 
$F_t(z)\in \ataz$, a \cNcs (noncommutative symmetric)
system (\cite{GTS-I}) $\Oft\in \cDazz[[t]]^{\times 5}$ 
has been constructed. 
Then, by the universal property of the 
\cNcs  system $(\cNsf, \Pi)$ over 
the Hopf algebra $\cNsf$ of NCSF's 
(noncommutative symmetric functions) 
(\cite{G-T}), we have a Hopf algebra 
homomorphism $\cSft: \cNsf \to \cDazz$.
Consequently, as pointed out in \cite{GTS-I}
as one of the main motivations 
for the introduction of 
the \cNcs  systems, by applying 
the homomorphism $\cSft$ to the identities 
of the NCSF in the \cNcs  system 
$\Pi$, we get a host of identities 
for the corresponding differential operators 
in the \cNcs  system $\Oft$.

In this paper, we first give a list of the identities 
for any two sequences of differential operators in the
\cNcs  system $\Oft$. These identities either  
come directly from the identities of the corresponding NCSF's 
derived in \cite{G-T} by applying the homomorphism $\cSft$, 
or can be derived easily from those identities 
of NCSF's by some simple arguments.
Secondly, by using these identities 
for the differential operators in $\Oft$ and 
the special forms of certain 
differential operators 
in $\Oft$ when $F_t(z)=z-tH(z)$ for some 
$H(z)\in \kzz^{\times n}$,
we derive some formulas in terms of differential operator 
in the system $\Oft$ for the Taylor series expansions of 
$u(F_t)$ and $u(F_t^{-1})$ $(u(z)\in \kttzz)$; 
the D-Log and the formal flow of $F_t$, 
and more importantly, 
some inversion formulas 
for the inverse maps of $F_t$. 
Finally, we discuss a connection 
of the well-known Jacobian conjecture 
with NCSF's.

Note that, the NCSF's were first introduced 
and studied in the seminal paper 
\cite{G-T} in $1994$.
\cNcs  systems over associative algebras 
were first formulated in \cite{GTS-I}, 
but mainly motivated 
by the introduction of the NCSF's 
in \cite{G-T} 
(see Definition \ref{Main-Def}).
Actually, in some sense, 
a \cNcs  system $\Omega$ over 
an associative $K$-algebra $A$ 
can be viewed as a system of analogs 
of the NCSF's in $A$ defined by 
Eqs.\,$(\ref{UE-0})$--$(\ref{UE-4})$ 
over $A$, which formally are same as 
the defining equations of certain NCSF's 
over the free $K$-algebra $\cNsf$ 
generated by a sequence of noncommutative 
free variables $\Lambda_m$ $(m\geq 1)$.
For some general discussions on 
the \cNcs  systems, see \cite{GTS-I}.
For more studies on NCSF's, 
see \cite{T}, \cite{NCSF-II}, 
\cite{NCSF-III}, \cite{NCSF-IV}, 
\cite{NCSF-V} and \cite{NCSF-VI}.

While, on the other hand, 
the inversion problem, 
which is mainly to study various 
properties of the inverse maps of analytic maps,
has much longer history. 
Since as early as $1770$ 
when L. Lagrange  \cite{L} proved 
the so-called Lagrange inversion formula,
there have been numerous papers 
devoted to find various 
inversion formulas, i.e. formulas
for inverse maps 
(see \cite{WZ}, \cite{BurgersEq}, \cite{NC-IVP} 
and references there).
The study on the inversion problem 
was greatly intensified since
O. H. Keller \cite{Ke} in $1939$  
proposed the well-known Jacobian 
conjecture which says, {\it  any polynomial 
map $F(z): K^{\times n} \to K^{\times n}$ 
with $j(F):=\det(\frac{\p F_i}{\p z_j})=1$ 
must be an automorphism of 
$K^{\times n}$ and its inverse map 
$G(z)=F^{-1}(z)$ 
must also be a polynomial map.} 
Despite enormous efforts 
from mathematicians in more than sixty years,
the conjecture is still open 
even for the case $n=2$. 
In 1998, S. Smale \cite{S} 
included the Jacobian conjecture 
in his list of $18$ important 
mathematical problems 
for $21$st century.
For more history and known results on
the Jacobian conjecture, 
see \cite{BCW}, \cite{E} and 
references there.
For some very recent developments 
on the conjecture, see \cite{BE1}, 
\cite{Me} and \cite{HNP}.

The arrangement of this paper is as follows.
In Section \ref{S2}, we mainly fix some notation 
and recall some results from \cite{GTS-II} 
that will be needed throughout this paper. 
In Subsection \ref{S2.1}, we briefly recall 
the \cNcs  systems in general and the universal 
\cNcs  system $(\cNsf, \Pi)$ formed 
by the generating functions of 
certain NCSF's  introduced in \cite{G-T}. 
In Subsection \ref{S2.2}, we recall 
the \cNcs systems $\Oft$ $(F_t\in \ataz)$
and the corresponding differential operator
specialization $\cSft\to \cDazz$ constructed 
in \cite{GTS-II}.  In Section \ref{S3}, 
we focus on the special automorphism 
$F_t\in \ataz$ such that $F_t(z)=z-tH(z)$ 
for some $H(z)\in \kzz^{\times n}$. 
We recall some results from \cite{NC-IVP} 
and \cite{GTS-II} which show that, 
in this case, certain differential operators 
have some simple forms.
Together with the specialization $\cSft$, 
the results in this section will be crucial 
for most of the formulas that will be  
derived in Subsections \ref{S5.3} 
and \ref{S5.4}, and also later a connection 
of NCSF's with the Jacobian conjecture 
in Subsection \ref{S5.5}.
In Section \ref{S4},
we mainly derive some
identities for the NCSF's in the
universal \cNcs system 
$(\cNsf, \Pi)$, which will be needed 
in next section. Note that, 
by the explicit correspondence 
in Corollary \ref{S-Correspondence-2}, 
applying the specialization 
$\cS_{F_t}: \cNsf\to \cDazz$ $(F_t\in \ataz)$
or simply changing the up case letters 
to the lower case letters, 
all the identities given in this section 
will become the identities of 
the corresponding differential 
operators in the \cNcs system 
$(\cDaz, \Oft)$.
In Subsection \ref{S4.1}, 
we give a list 
of the identities between any two 
sequences of the NCSF's in 
the universal \cNcs  system
$(\cNsf, \Pi)$. In Subsection \ref{S4.2}, 
we let $u$ be another formal central parameter
and derive the formulas for 
$e^{-u \Phi(t)}$ in terms of the NCSF's in $\Pi$.  
In Section \ref{S5}, 
we mainly apply the identities of NCSF's derived 
in the previous section and
the specialization 
$\cS_{F_t}$ $(F_t\in \ataz)$ 
in Theorem \ref{S-Correspondence}
to derive some formulas
for the inverse maps, the D-Log's and 
the formal flow of the automorphism 
$F_t\in \ataz$. 
In Subsections \ref{S5.1} and \ref{S5.2}, 
we derive formulas for the D-Log of and 
the formal flow generated by $F_t\in \ataz$, 
respectively, in terms of the differential 
operators in the \cNcs system $\Oft$.
In Subsection \ref{S5.3},
we mainly apply the results in the previous 
two subsections to the special automorphisms 
$F_t(z)\in \etaz$ to 
derive some inversion formulas. 
In Subsection \ref{S5.4}, motivated 
by the fact that 
$C_m(z) \in \kzz^{\times n}$
$(m\geq 1)$ in the commutative case
capture the nilpotence 
of the Jacobian matrix $JH$, 
we give formulas for 
$C_m(z)$'s in terms of 
the differential operators in $\Oft$ 
other than $\psi_m$'s.
Finally, in Subsection \ref{S5.5}, 
we discuss a connection of NCSF's with 
the well-known Jacobian conjecture.

\renewcommand{\theequation}{\thesection.\arabic{equation}}
\renewcommand{\therema}{\thesection.\arabic{rema}}
\setcounter{equation}{0}
\setcounter{rema}{0}

\section{\bf Differential Operator Specializations of NCSF's} \label{S2}

Let $K$ be any unital commutative $\bQ$-algebra and 
$A$ any unital associative but not necessarily commutative 
$K$-algebra. Let $t$ be a formal central parameter, 
i.e. it commutes with all elements of $A$, and $A[[t]]$ 
the $K$-algebra of formal power series 
in $t$ with coefficients in $A$. 
First let us recall the following notion formulated 
in \cite{GTS-I}.

\begin{defi} \label{Main-Def}
For any unital associative  $K$-algebra $A$, a $5$-tuple $\Omega=$ 
$( f(t)$, $g(t)$, $d\,(t)$, $h(t)$, $m(t) ) 
\in A[[t]]^{\times 5}$ is said 
to be a {\it \cNcs $($Noncommutative Symmetric$)$ system}
over $A$ 
if the following equations are satisfied.
\allowdisplaybreaks{
\begin{align}
&f(0)=1 \label{UE-0}\\
& f(-t)  g(t)=g(t)f (-t)=1, \label{UE-1}   \\
& e^{d\,(t)} = g(t), \label{UE-2} \\
& \frac {d g(t)} {d t}= g(t) h(t), \label{UE-3}\\ 
& \frac {d g(t)}{d t} =  m(t) g(t).\label{UE-4}
\end{align} }
\end{defi}

When the base algebra $K$ is clear in the context, we also call 
the ordered pair $(A, \Omega)$ a {\it \cNcs  system}. 
Since \cNcs  systems often come from generating functions 
of certain elements of $A$ that are under concern, 
the components of $\Omega$ will also be refereed as 
the {\it generating functions} of their coefficients. 

In this section, we mainly fix necessary notations and 
recall some results from \cite{GTS-I} 
and \cite{GTS-II} that will be needed later.
In Subsection \ref{S2.1}, we briefly recall 
the \cNcs  system $(\cNsf, \Pi)$ formed by 
generating functions of some of the NCSF's 
defined in \cite{G-T} and 
its universal property (see Theorem \ref{Universal}).  
In Subsection \ref{S2.2}, we recall
the \cNcs systems (\cite{GTS-II}) 
over differential operator algebras
and the resulted correspondence
between NCSF's and the differential operators 
in the system.

\subsection{The Universal \cNcs  System from
Noncommutative Symmetric Functions}\label{S2.1}

Let $\Lambda=\{ \Lambda_m\,|\, m\geq 1\}$ 
be a sequence of noncommutative 
free variables and $\cNsf$ 
the free associative algebra 
generated by 
$\Lambda$ over $K$.  For convenience, 
we also set $\Lambda_0=1$.
We denote by
$\lambda (t)$ the generating function of 
$\Lambda_m$ $(m\geq 0)$, i.e. we set
\begin{align}\label{lambda(t)}
\lambda (t):= \sum_{m\geq 0} t^m \Lambda_m 
=1+\sum_{k\geq 1} t^m \Lambda_m.
\end{align}

In the theory of NCSF's (\cite{G-T}), 
$\Lambda_m$ $(m\geq 0)$ is 
the noncommutative analog 
of the $m^{th}$ classical (commutative) 
elementary symmetric function 
and is called the {\it $m^{th}$ 
$(\text{noncommutative})$ 
elementary symmetric function.}

To define some other NCSF's, we consider 
Eqs.\,$(\ref{UE-1})$--$(\ref{UE-4})$ 
over the free $K$-algebra $\cNsf$
with $f(t)=\lambda(t)$. The 
solutions for $g(t)$, $d\,(t)$, 
$h(t)$, $m(t)$ exist and are unique, 
whose coefficients will be the NCSF's 
that we are going to define.
Following the notation in \cite{G-T} 
and \cite{GTS-I}, we denote the resulted 
$5$-tuple by 
\begin{align}
\Pi=(\lambda(t),\, \sigma(t),\, \Phi(t),\, \psi(t),\, \xi(t))
\end{align}
and write the last 
four generating functions of 
$\Pi$ explicitly as follows.

\allowdisplaybreaks{
\begin{align}
\sigma (t)&=\sum_{m\geq 0} t^m S_m,  \label{sigma(t)} \\
\Phi (t)&=\sum_{m\geq 1} t^m \frac{\Phi_m}m  \label{Phi(t)}\\
\psi (t)&=\sum_{m\geq 1} t^{m-1} \Psi_m, \label{psi(t)}\\
\xi (t)&=\sum_{m\geq 1} t^{m-1} \Xi_m.\label{xi(t)}
\end{align}}

Following \cite{G-T},
we call $S_m$ ($m\geq 1$) the 
{\it $m^{th}$ $(\text{noncommutative})$ complete 
homogeneous symmetric function} and
$\Phi_m $ (resp.\,\,$\Psi_m$) 
the {\it $m^{th}$ power sum symmetric function 
of the second $($resp.\,\,first$)$ kind}. 
Following \cite{GTS-I}, 
we call $\Xi_m \in \cNsf$ $(m\geq 1)$ 
the {\it $m^{th}$ $(\text{noncommutative})$ 
power sum symmetric function of the third kind}.

The following  proposition proved in \cite{G-T} 
and \cite{NCSF-II} will be very useful 
for our later arguments.


\begin{propo}\label{omega-Lambda}
Let $\omega_\Lambda$ be the anti-involution of 
$\cNsf$ 
which fixes $\Lambda_m$ $(m\geq 1)$.
Then, for any $m\geq 1$, we have
\begin{align}
\omega_\Lambda (S_m)&=S_m,  \label{omega-Lambda-e1}\\
\omega_\Lambda (\Phi_m)&=\Phi_m,\label{omega-Lambda-e2} \\
\omega_\Lambda (\Psi_m)&=\Xi_m. \label{omega-Lambda-e3}
\end{align}
\end{propo}

By applying Proposition $2.10$ in \cite{GTS-I} to the 
universal \cNcs system $(\cNsf, \Pi)$, 
we have the following proposition.

\begin{propo}\label{tau-flip}
Let $\tau$ be the involution of $\cNsf$ 
such that $\tau (\Phi_m)=-\Phi_m$ for any $m\geq 1$.
Then, we have  
\begin{align}
\tau (\Lambda_m)&=(-1)^m S_m   \, ,\label{tau-flip-e1} \\
\tau (S_m)&=(-1)^m \Lambda_m \, ,\label{tau-flip-e2} \\
\tau (\Psi_m)&=- \Xi_m \, ,\label{tau-flip-e3}  \\
\tau (\Xi_m)&=-\Psi_m\, .\label{tau-flip-e4} 
\end{align}
\end{propo}

Next, let us recall the following 
$K$-Hopf algebra structure 
of $\cNsf$. It has been shown in 
\cite{G-T} that $\cNsf$ is the universal enveloping algebra 
of the free Lie algebra generated 
by $\Psi_m$ $(m\geq 1)$. Hence, it has a Hopf  
$K$-algebra structure as all other 
universal enveloping algebras 
of Lie algebras do. Its co-unit $\epsilon:\cNsf \to K$,
 co-product $\Delta$ and 
 antipode $S$ are uniquely determined by 
\begin{align}
\epsilon (\Psi_m)&=0, \label{counit} \\
\Delta (\Psi_m) &=1\otimes \Psi_m 
+\Psi_m\otimes 1, \label{coprod}\\
S(\Psi_m) & =-\Psi_m,\label{antipode}
\end{align}
for any $m\geq 1$. 

From the definitions of the NCSF's above, 
we see that $(\cNsf, \Pi)$ obviously forms a \cNcs  system.
More importantly, as shown in Theorem $2.1$ in \cite{GTS-I}, 
we have the following important theorem on 
the \cNcs  system $(\cNsf, \Pi)$. 

\begin{theo}\label{Universal}
Let $A$ be a $K$-algebra and $\Omega$ 
a \cNcs  system over $A$. Then, 
There exists a unique $K$-algebra homomorphism 
$\cS: \cNsf\to A$ such that 
$\cS^{\times 5} (\Pi)=\Omega$.
\end{theo}

Note that, when $A$ is further a $K$-bialgebra $($resp.\,\,$K$-Hopf algebra$)$ 
some sufficient conditions for the homomorphism 
$\cS: \cNsf\to A$ in the theorem above to be  
a homomorphism of $K$-bialgebras 
$($resp.\,\,$K$-Hopf algebras$)$ 
were also given in Theorem $2.1$ in \cite{GTS-I}. 

\begin{rmk}\label{comm-case}
By taking the quotient over 
the two-sided ideal generated by the commutators 
of $\Lambda_m$'s or applying a similar argument for 
proof of Theorem \ref{Universal}, 
it is easy to see that, 
over the category of commutative $K$-algebras, 
the universal \cNcs  system 
is given by the generating functions of
the corresponding classical 
$($commutative$)$ symmetric 
functions \cite{Mac}.
\end{rmk}

\subsection{\cNcs  Systems over Differential Operator Algebras}
\label{S2.2}

In this subsection, we briefly recall 
the \cNcs  systems constructed in \cite{GTS-II} over 
the differential operator algebras in commutative 
or noncommutative free variables.
First, let us fix the following notation.

Let $K$ be any unital 
commutative $\bQ$-algebra as before
and $z=(z_1, z_2, ... , z_n)$ commutative 
or noncommutative free 
variables.\footnote{Since most of the results 
as well as 
their proofs in this paper 
do not depend on the commutativity 
of the free variables $z$,  
we will not distinguish the commutative 
and the noncommutative case, 
unless stated otherwise, 
and adapt the notations 
for noncommutative variables 
uniformly for the both cases.}
Let $t$ be a formal central parameter, 
i.e. it commutes with $z$ and elements of $K$.
We denote by $\kzz$ and 
$\kttzz$ the $K$-algebras of 
formal power series in $z$ over  
$K$ and $K[[t]]$, respectively.

By a {\it $K$-derivation} or simply {\it derivation} 
of $\kzz$, we mean a $K$-linear $\delta: \kzz\to \kzz$ 
that satisfies the Leibniz rule,
i.e. for any $f, g\in \kzz$, we have
 \begin{align}\label{Leibniz}
\delta (fg)=(\delta f)g+f(\delta g).
\end{align}
We will denote by 
$\cDrkzz$ or $\cDrzz$, when the base algebra
$K$ is clear from the context, 
the set of all $K$-derivations 
of $\kzz$.
The unital subalgebra of 
$\text{End}_k(\kzz)$
generated by all
$K$-derivations 
of $\kzz$ will be denoted by 
$\cDkzz$ or $\cDzz$. 
Elements of $\cDkzz$ will be called 
{\it $(\text{formal})$ 
differential operators} in the
commutative and 
noncommutative variables $z$.

For any $\alpha\geq 1 $, 
we denote by $\cDrazz$ 
the set of the $K$-derivations 
of $\kzz$ which increase 
the degree in $z$ by 
at least $\alpha-1$. 
The unital subalgebra of 
$\cDzz$ generated by elements of 
$\cDrazz$ will be denoted by $\cDazz$. 
Note that, by the definitions above,
the operators of scalar multiplications
are also in $\cDzz$ and $\cDazz$.
When the base algebra is $K[[t]]$ 
instead of $K$ itself,
the notation 
$\cDrzz$, $\cDzz$, $\cDrazz$ 
and $\cDazz$ will be denoted by  
$\cDrtzz$, $\cDtzz$, $\cDrtazz$ 
and $\cDtazz$, respectively.
For example,  $\cDrtazz$ stands for 
the set of all $K[[t]]$-derivations of $\kttzz$ 
which increase the degree in $z$ 
by at least $\alpha-1$. 
Note that, $\cDrtazz=\cDrazz[[t]]$ and 
$\cDtazz=\cDazz[[t]]$.

For any $1\leq i\leq n$ and $u(z)\in \kzz$, 
we denote by $\lb u(z) \fr \p{\p z_i}\rb $
the $K$-derivation which maps $z_i$ to $u(z)$ and $z_j$ to $0$ 
for any $j\neq i$. 
For any $\Vec{u}=(u_1, u_2, \cdots, u_n)\in \kzz^{\times n}$, 
we set 
\begin{align}\label{Upz}
[\Vec{u}\pz]:=\sum_{i=1}^n [u_i \fr\p{\p z_i}].  
\end{align}

Note that, in the noncommutative case, 
we in general do {\bf not} have
$\lb u(z) \fr \p{\p z_i}\rb  g(z)  = u(z)  
\fr {\p g}{\p z_i}$ for all 
$u(z), g(z)\in \kzz$. This is the reason why
we put a bracket $[\cdot]$ in the notation above 
for the $K$-derivations.
With this notation, it is easy to see that any 
$K$-derivations $\delta$ 
of $\kzz$ can be written uniquely  as 
$\sum_{i=1}^n \lb f_i(z)\fr\p{\p z_i}\rb$ 
with $f_i(z)=\delta z_i\in \kzz$ 
$(1\leq i\leq n)$.

With the commutator bracket, 
$\cDrazz$ $(\alpha\geq 1)$ forms a Lie algebra 
and its universal enveloping algebra is exactly 
the differential operator algebra
$\cDazz$.  Consequently, $\cDazz$ $(\alpha\geq 1)$
has a Hopf algebra structure as all 
other enveloping algebras of Lie algebras do.
In particular,
Its coproduct $\Delta$, antipode $S$ and co-unit $\epsilon$
are uniquely determined by the properties
\begin{align}
\Delta(\delta)&= 1\otimes\delta+\delta\otimes 1,\label{Coprd-delta} \\
S(\delta)&=-\delta, \label{antipd-delta}\\
\epsilon (\delta) &=\delta \cdot 1, \label{Counit-delta}
\end{align}
respectively, for any $\delta \in \cDrzz$.

For any $\alpha\geq 1 $, let $\ataz$ 
be the set of all the automorphism 
$F_t(z)$ of $\kttzz$ over $K[[t]]$, 
which have the form
$F(z)=z-H_t(z)$ for some 
$H_t(z)\in \kttzz^{\times n}$ 
with $o(H_t(z))\geq \alpha$ and $H_{t=0}(z)=0$. 
Note that, for any $F_t \in \ataz$ as above, 
its inverse map $G_t:=F_t^{-1}$
can always be written uniquely 
as $G_t(z)=z+M_t(z)$ \label{Mt(z)}
for some $M_t(z)\in \kttzz^{\times n}$ 
with $o(M_t(z))\geq \alpha$ 
and $M_{t=0}(z)=0$. 

Now we recall the \cNcs  systems constructed 
in \cite{GTS-II} over 
the differential operator algebras 
$\cDazz$ $(\alpha\geq 1)$. 
We fix an $\alpha \geq 1$ and 
an arbitrary $F_t\in \ataz$. 
We will always 
let $H_t(z)$, $G_t(z)$ 
and $M_t(z)$ be determined 
as above. 
The \cNcs  system
\begin{align}\label{Def-Omega-Ft}
\Omega_{F_t}=(f(t)\, g(t),\, d(t),\, 
h(t),\, m(t)) \in \cDazz[[t]]^{\times 5}.
\end{align}
is determined as follows.

The last two components are given directly by the following two 
$K[[t]]$-derivations of the $K[[t]]$-algebra $\kttzz$.

\begin{align}
h(t)&:=\left[\frac{\p M_t}{\p t}(F_t) \pz \right ],\label{Def-h(t)-0}\\
m(t)&:=\left[\frac{\p H_t}{\p t}(G_t) \pz \right ]. \label{Def-m(t)-0}
\end{align}

The first three components are given 
by the following proposition which was
proved in Section $3.2$ in \cite{GTS-II}. 

\begin{propo}\label{TaylorExpansion-DLog}
There exist unique $f(t), g(t), d(t) \in \cDtazz$ 
with $f(0)=1$ and $d(0)=0$
such that, for any $u_t(z)\in \kttzz$, we have
\begin{align}
f(-t) u_t(z)&=u_t(F_t),\label{NewTaylorExpansion-e1}\\
g(t) u_t(z)&=u_t(G_t), \label{NewTaylorExpansion-e2}\\
e^{d(t)} u_t(z) &=u_t(G_t), \label{NewDLog-e2}
\end{align}
where, as usual, the exponential in 
Eq.\,$(\ref{NewDLog-e2})$ is given by 
\begin{align}\label{L2.2.3-e2}
e^{d(t)} = \sum_{m\geq 0} \frac {d(t)^m}{m!}.
\end{align}
\end{propo}

By using the Taylor series expansions of $u(F_t)$
and $u(G_t)$, it is easy to show 
that the differential operators 
$f(t)$ and $g(t)$ can be given as follows.
 
\begin{lemma}\label{Taylor-fg(t)}
\begin{align}
f(t)&= \sum_{m \geq 0} \frac {1}{m!}
\left.  \lb H_t(w) \pz \rb^m \,  \right |_{w=z}, \label{Taylor-fg(t)-e1} \\
g(t)&= \sum_{m \geq 0} \frac {1}{m!}
\left.  \lb M_t(w) \pz \rb^m \, \right |_{w=z},  \label{Taylor-fg(t)-e2}
\end{align}
where $w=(w_1, \cdots, w_n)$ are free variables that are 
independent with $z$,  and the notation $|_{w=z}$ 
simply means that, after applying 
the differential operator before $|_{w=z}$ 
and then replacing $w$ back by $z$.
\end{lemma}

Note that, when we write $d(t)$ above as 
$d(t)= -\lb a_t(z)\pz \rb$ for some 
$a_t(z)\in t\kttzz$, then we get the so-called 
{\it D-Log}\label{D-Log} $a_t(z)$ of the automorphism 
$F_t(z)\in \ataz$, which has been studied 
in \cite{E1}--\cite{E3}, \cite{Z-exp} and \cite{WZ} 
for the commutative case. 

We define five sequences of 
 differential operators by writing the 
 components of $\Omega_{F_t}$ as follows.
 \allowdisplaybreaks{
 \begin{align}
f(t)&= \sum_{m=0}^\infty t^m \lambda_m, \label{Def-f(t)} \\
g(t)&= \sum_{m=0}^\infty t^m s_m. \label{Def-g(t)} \\
d(t)&=\sum_{m=1}^\infty \frac {t^{m}}m \phi_m \label{Def-d(t)} \\
h(t)&=\sum_{m\geq 1} \psi_m t^{m-1},\label{Def-h(t)} \\
m(t)&=\sum_{m\geq 1} \xi_m t^{m-1}. \label{Def-m(t)}
\end{align}}

Throughout this paper, 
we will also denote each sequence of 
the differential operators above 
by the corresponding letter without sub-index. 
For example, $\lambda$ denotes the sequence of 
the differential operator 
$\{\lambda_m\,|\, m\geq 0\}$ defined in 
Eq.\,(\ref{Def-f(t)}).

Let $\cSft: \cNsf\to \cDazz$ be the unique $K$-algebra 
homomorphism that maps $\Lambda_m$ to $\lambda_m$ 
for any $m\geq 1$. Note that, since $\cNsf$ is 
freely generated by $\lambda_m$ $(m\geq 1)$, 
the map $\cSft$ is well-defined. 
The main result we need later 
is the following theorem proved in \cite{GTS-II}.

\begin{theo} \label{S-Correspondence} 
For any $\alpha\geq 1 $ and $F_t(z)\in \ataz$, we have,

$(a)$ the $5$-tuple $\Omega_{F_t}$ defined 
in Eq.\,$(\ref{Def-Omega-Ft})$
forms a \cNcs system over 
the differential operator algebra $\cDazz$.

$(b)$ $\cS_{F_t}: \cNsf \to \cDazz$ defined above 
is the unique $K$-Hopf algebras homomorphism 
such that $\cS_{F_t}^{\times 5}(\Pi)=\Omega_{F_t}$. 
\end{theo}

\begin{corol} \label{S-Correspondence-2}
For any $m\geq 1$, we have the following differential 
operator realizations for the corresponding NCSF's.
 \begin{align}
 \mathcal S_{F_t} (\Lambda_m) & =\lambda_m, \label{L-l} \\
 \mathcal S_{F_t} (S_m) & =s_m, \label{S-s} \\
 \mathcal S_{F_t} (\Psi_m) & =\psi_m, \label{Psi-psi}  \\
 \mathcal S_{F_t} (\Phi_m)& =\phi_m, \label{Phi-phi}  \\
 \mathcal S_{F_t} (\Xi_m)& =\xi_m. \label{Xi-xi} 
 \end{align}
 \end{corol}

For more properties of the specialization 
$\cSft:\cNsf\to \cDazz$, see \cite{GTS-II} 
and \cite{GTS-V}.
Finally, let us point out the following 
result in \cite{GTS-V} that will be helpful 
to our later arguments.

For any $z$ and $\alpha\geq 1$ as before, 
we let $\mathbb B^{[\alpha]}_t\langle z \rangle$
be the set of automorphisms $F_t=z-H_t(z)$ 
of the polynomial algebra 
$\ktz$ over $K[t]$ such that 
the following conditions are satisfied.
\begin{enumerate}
\item[$\bullet$]  $H_{t=0}(z)=0$. 
\item[$\bullet$] 
$H_t(z)$ is homogeneous in $z$ 
of degree $d \geq \alpha$.
\item[$\bullet$] 
With a proper permutation of the free variables $z_i$'s,
the Jacobian matrix $JH_t(z)$ becomes strictly lower triangular.
\end{enumerate}

\begin{theo}\label{StabInjc-best}
In both commutative and noncommutative cases, 
the following statement holds. 
\vskip2mm
For any fixed $\alpha\geq 1$ and non-zero $P \in \cNsf$, 
there exist $n\geq 1$ 
$($the number  of the free variable 
$z_i$'s$)$ and $F_t(z)\in \mathbb B^{[\alpha]}_t\langle z\rangle$ 
such that $\cS_{F_t} (P)\neq 0$.
\end{theo}

\renewcommand{\theequation}{\thesection.\arabic{equation}}
\renewcommand{\therema}{\thesection.\arabic{rema}}
\setcounter{equation}{0}
\setcounter{rema}{0}

\section{\bf A Family of Special Automorphisms $F_t(z)$} \label{S3}

Let $K$, $z$, $t$ and $\ataz$ $(\alpha\geq 1)$
as fixed in Subsection \ref{S2.2}.
We will also freely use other notations fixed 
in the earlier sections. 
First, let us set $\etaz$ be the set of 
all automorphisms $F_t\in \ataz$
such that $F_t(z)=z-tH(z)$
for some $H(z)\in \kzz^{\times n}$.

Note that, the automorphisms 
$F_t(z)\in \etaz$ can be viewed as 
special deformations of 
the automorphisms $F(z)=z-H(z)$ 
of the $K$-algebra $\kzz$. They
have been studied in \cite{BurgersEq} 
for the commutative case and later 
in \cite{NC-IVP} for
the noncommutative case.

In this section, we mainly recall 
some results from \cite{GTS-I}, 
which show that, for the automorphisms 
$F_t(z)\in \etaz$, the differential 
operators $\lambda_m$'s, $\psi_m$'s 
and $\xi_m$'s in the 
\cNcs  system $\Oft$ have 
some special forms. Together with the correspondence
in Corollary \ref{S-Correspondence-2}, 
the results in this section are responsible 
for most of the formulas that will be derived in  
Subsections \ref{S5.3} and  \ref{S5.4}, 
and also later a connection of NCSF's 
with the Jacobian conjecture 
in Subsection \ref{S5.5}.

First, let us fix a $F_t\in \etaz$ and always write
the formal inverse map $G_t(z):=F^{-1}_t(z)$ as 
$G_t(z)=z+tN_t(z)$ with $N_t(z)\in \ktzz^{\times n}$.
Note that, in terms of the notation 
in Section \ref{S2.2}, 
we have 
\begin{align}
H_t(z)&=tH(z), \label{Special-H} \\
M_t(z)&=tN_t(z). \label{Special-N}
\end{align}

Furthermore,  we define a sequence 
$\{N_{[m]}(z) \in \kzz^{\times n} | \, m\geq 0\}$ 
by writing 
\begin{align}\label{Def-Nm}
N_t(z)= \sum_{m=0}^{+\infty} t^{m-1} N_{[m]}(z).
\end{align}

By Lemma $4.1$ in \cite{NC-IVP}, the differential 
operators defined in Eqs.\,(\ref{Def-m(t)}) 
and (\ref{Def-h(t)}) have the following 
special forms. 

\begin{lemma} \label{L4.1.1}
\begin{align}
m(t)&=\lb N_t(z)\pz \rb, \label{L4.1.1-e1} \\
h(t)&=\sum_{m\geq 1}  t^{m-1} \lb C_m(z)\pz \rb, \label{L4.1.1-e2}
\end{align}
where $C_m(z)\in \kzz^{\times n}$ $(m\geq 1)$ 
are defined recursively by
\begin{align}
C_1(z)& =H(z),\label{L4.1.1-e3}\\
C_m(z)&=\lb C_{m-1}(z)\pz \rb H,\label{L4.1.1-e4}
\end{align}
for any $m\geq 2$.
\end{lemma}
 
Note that, by the mathematical induction, 
it is easy to show that, 
when $z$ are commutative variables, 
we further have
\begin{align}
C_m(z)=(JH)^{m-1}H \label{Special-Cm}
\end{align}
for any $m\geq 1$.

Consequently,  the $K$-derivations 
$\psi_m$ and $\xi_m$ $(m\geq 1)$ defined in 
Eqs.\,(\ref{Def-m(t)}) and (\ref{Def-h(t)})
have the following simple forms.

\begin{corol}\label{C4.1.2}
For any $m\geq 1$, we have
\begin{align}
\psi_m &=\lb C_m(z)\pz \rb, \label{Special-psi-1}\\
\xi_m  &=\lb N_{[m]}(z)\pz \rb. \label{Special-xi}
\end{align}

When $z$ are commutative variables, 
we further have
\begin{align}
\psi_m =\lb \lp (JH)^{m-1} H \rp \pz \rb \label{Special-psi-m}
\end{align}
\end{corol}

Note that, by applying Eq.\,(\ref{Taylor-fg(t)-e1})
to $F_t\in \etaz$, we have the following special form
for the differential operators $\lambda_m$'s.

\begin{lemma} \label{Special-lambda}
For any $m\geq 1$ and $u(z)\in \kzz$, we have
\begin{align}\label{Special-lambda-1}
\lambda_m u(z) = \frac 1{m!} \left. \lb  H(w) \pz \rb^m u(z)\,  \right |_{w=z} 
\end{align}
In a short notation as in Lemma \ref{Taylor-fg(t)}, we have, 
\begin{align}\label{Special-lambda-2}
\lambda_m  &= \frac 1{m!} \left. \lb  H(w) \pz \rb^m \,  \right |_{w=z}\, . 
\end{align}
In particular, when $z$ are commutative free variables, we have
\begin{align}\label{Special-lambda-3}
\lambda_m  =  \sum_{\substack{I\in \bN^n \\ |I|=m }} 
\frac 1{I!} H^I(z) \frac{\p^m} {\p z^I}.
\end{align}
\end{lemma}

\renewcommand{\theequation}{\thesection.\arabic{equation}}
\renewcommand{\therema}{\thesection.\arabic{rema}}
\setcounter{equation}{0}
\setcounter{rema}{0}

\section{\bf Identities of Noncommutative Symmetric Functions} \label{S4}

In this section, 
we mainly derive some
identities for the NCSF's in the
universal \cNcs system 
$(\cNsf, \Pi)$, which will be needed 
in next section. Note that, 
by the explicit correspondence 
in Corollary \ref{S-Correspondence-2}, 
applying the specialization 
$\cS_{F_t}: \cNsf\to \cDazz$ $(F_t\in \ataz)$
or simply changing the up case letters 
to the lower case letters, 
all the identities given in this section 
will become the identities of 
the corresponding differential 
operators in the \cNcs system 
$(\cDaz, \Oft)$.
In Subsection \ref{S4.1}, 
we give a list 
of the identities between any two 
sequences of the NCSF's in 
the universal \cNcs  system
$(\cNsf, \Pi)$. In Subsection \ref{S4.2}, 
we let $u$ be another formal central parameter
and derive the formulas for 
$e^{-u \Phi(t)}$ in terms of the NCSF's in $\Pi$.  
These formulas will be used in 
Subsection \ref{S5.2} 
to derive certain formulas for 
the formal flow generated by 
$F_t\in \ataz$.

First, let us fix the 
following notations.

\vskip2mm
{\bf Notation:}
\vskip2mm
\begin{enumerate}

\item[(i)] For any composition $I$, i.e. 
an ordered finite sequence 
$I=(i_1, i_2, ... , i_m)$ of 
positive integers, we define the {\it length} 
$\ell(I)$ of $I$ to be $m$ and the {\it weight} 
$|I|$ of $I$ to be $\sum_{j=1}^m i_j$.
We denote by $\mathcal C$ (resp.\,\,$\cC_m$ $(m\geq 1)$) 
the set of all compositions $I$ (resp.\,\,with $|I|=m$).

\item[(ii)] Let $I$ as given in (i) and
$\{w_m\,|\, m\geq 1\}$ a sequence of letters or 
elements of a $K$-algebra, 
we set $w^I=w_{i_1}w_{i_2}\cdots w_{i_m}$. 

\item[(iii)] Let $I$ as given in $(i)$, we define
the {\it first part}\, $lp(I)$ 
(resp.\,\,the {\it last part}\, $lp(I)$) 
of $I$ to be $i_1$ (resp.\,\,$i_m$).
Furthermore, we also set 
\allowdisplaybreaks{
\begin{align}
\pi(I)&=\prod_{j=1}^m i_j, \\
\pi_u(I)&=\prod_{j=1}^m (i_1+\cdots +i_j),\\
 sp (I)&=\ell (I)!\, \pi(I).
\end{align}}

\item[(iv)] For any fixed composition as in (i), 
we define  the {\it mirror image} $\bar I$ of $I$ to 
be the composition obtained by 
reversing the ordered sequence $I$, 
i.e. $\bar{I}:=(i_m, i_{m-1}, \cdots, i_1)$. 

\item[(v)] For any compositions $I, J\in \cC$, we denote by $I\cdot J$ 
the {\it concatenation product} of $I$ and $J$. For example, 
if $I=(3, 2, 4)$ and $J=(5, 7)$, then 
$I\cdot J=(3, 2, 4, 5, 7)$.

\item[(vi)] Let $I$ be as in (i) and 
$J=(j_1, j_2, \cdots, j_k)$ another composition.
We say $J$ is a {\it refinement} of $I$, 
denoted by  $J\succcurlyeq I$ or $I\preccurlyeq J$,
if there exist $d_0=1< d_1< d_2 < \cdots < d_m=k+1$ 
such that, for any $1 \leq a \leq m$, 
we have 
$$
j_{d_{a-1}} + j_{d_{a-1}+1}+ \cdots + j_{d_{a}-1} = i_a.
$$
For example,  $(4,2,5,2,1) \succcurlyeq (6,5,3)$ and 
$(4,1,7)\preccurlyeq (3, 1, 1, 5, 2)$.

\item[(vii)] Let $I$ and $J$ be any two composition 
with $J\succcurlyeq I$. With notation fixed in (vi), we set , 
for any $1 \leq a \leq m$, 
$$
J_a:=\{\, j_{d_{a-1}}, \, j_{d_{a-1}+1}, \, \cdots ,\, j_{d_{a}-1}\, \}.
$$
Then we further set
\allowdisplaybreaks{
\begin{align}
 \ell(J, I)&:=\prod_{a=1}^m  \ell(J_a), \label{l-JI} \\
\pi_u(J, I)&:=\prod_{a=1}^m \pi_u (J_a),\label{pi-u-JI}\\ 
sp (J, I)&=\prod_{i=1}^m sp (J_i), \\ 
lp (J, I)&:=\prod_{a=1}^m lp (J_a), \label{lp-JI}\\ 
fp (J, I)&:=\prod_{a=1}^m fp (J_a). \label{fp-JI}
\end{align}}
\end{enumerate}

\subsection{Identities of the NCSF's in the \cNcs  System $\Pi$ }
\label{S4.1}

In this subsection, we give a list of the identities 
of the NCSF's in the \cNcs  System $\Pi$. 
Note that,  by simply applying the specialization 
$\cSft: \cNsf \to\cDazz$ in 
Theorem \ref{S-Correspondence} or just replacing 
the up case letters by the lower case letters, 
these identities will become the identities of
the differential operators in the \cNcs  system 
$\Omega_{F_t}$. 

First, we fix a composition $I=(i_1, i_2, ... , i_m)$ and 
start with the following five pairs of the identities of NCSF's,  
which have been derived in $\S 4$ of \cite{G-T}.

{\it
\begin{enumerate}

\item[$\bullet$] The relations between $\Lambda$ and $S$:  
\begin{align}
S^I & = \sum_{J\succeq I} (-1)^{ \ell (J)-|I| } \,\, \Lambda^J\, , \label{Lambda-S} \\
\Lambda^I & = \sum_{J \succeq I} (-1)^{\ell(J)-|I|}\, \, S^J\, .\label{S-Lambda}
\end{align}

\item[$\bullet$] The relations between $\Lambda$ and $\Psi$: 
\begin{align}
\Lambda^I & = (-1)^{|I|} \sum_{J\succeq I}  
\frac{(-1)^{\ell (J)}}{\pi_u( {\bar J}, {\bar I})}\, \, \Psi^J \, , \label{Psi-Lambda} \\
\Psi^I & = (-1)^{|I|} \sum_{J \succeq I} 
(-1)^{\ell(J)} fp (J, I) \, \, \Lambda^J\, . \label{Lambda-Psi}
\end{align}

\item[$\bullet$] The relations between $S$ and $\Psi$:
\begin{align}
S^I & =  \sum_{J\succeq I}  \frac{1}{\pi_u(  J, I)}\, \,  \Psi^J\, , \label{Psi-S}  \\
\Psi^I & = (-1)^{\ell(I) }\sum_{J \succeq I} (-1)^{ \ell(J)} lp (J, I) \, \, S^J\, . \label{S-Psi}
\end{align}

\item[$\bullet$] The relations between $\Lambda$ and $\Phi$: 
\begin{align}
\Lambda^I & =  (-1)^{|I|}\sum_{J\succeq I} 
 \frac{ (-1)^{ \ell(J)}  } {sp (  J, I)}\, \,  \Phi^J\, , \label{Phi-Lambda} \\
\Phi^I & =  (-1)^{|I|}\sum_{J\succeq I}  
 \frac{ (-1)^{ \ell(J)}\pi(I) }{\ell (  J, I)}\, \,  \Lambda^J\, . \label{Lambda-Phi}
\end{align}

\item[$\bullet$] The relations between $S$ and $\Phi$: 
\begin{align}
S^I & =  \sum_{J\succeq I}  \frac{1}{sp ( J, I)}\,\,  \Phi^J\, , \label{Phi-S}  \\
\Phi^I & = (-1)^{\ell(I) }
\sum_{J \succeq I}  \frac{(-1)^{ \ell(J)} \pi(I)}{\ell(J, I)}\,\,  S^J\, . \label{S-Phi}
\end{align}

\item[$\bullet$] The relations between $\Psi$ and $\Phi$: 
\begin{align}
\Phi^I & =(-1)^{\ell(I)}
\sum_{J \succeq I} \lp 
\sum_{J \succeq K \succeq I} 
\frac{(-1)^{ \ell(K)} \pi(I)}{\pi_u (J, K)\, \ell(K, I)} \rp \,   \Psi^J\, , \label{Psi-Phi}  \\
\Psi^I & = (-1)^{\ell(I) } \sum_{J \succeq I}
\lp \sum_{J \succeq K \succeq I}   \frac {(-1)^{\ell(K)}\, lp (K, I)}{sp (J, K)} \rp \,  \Phi^J\, .
\label{Phi-Psi}
\end{align}
\end{enumerate} }

The last two identities were not given explicitly in \cite{G-T}, 
but can be easily derived as follows.

\pf First, by combining Eq.\,(\ref{S-Phi}) with Eq.\,(\ref{Psi-S}), 
we have
\begin{align*}
\Phi^I & = (-1)^{\ell(I) }
\sum_{K \succeq I} (-1)^{ \ell(K)} \frac{\pi(I)}{\ell(K, I)}\, S^K \\
& =(-1)^{\ell(I) }
 \sum_{K \succeq I} (-1)^{ \ell(K)} \frac{\pi(I)}{\ell(K, I)} 
\sum_{J \succeq K}  \frac{1}{\pi_u (J, K)}\,  \Psi^J \\
& =(-1)^{\ell(I) }
\sum_{J \succeq K} \lp 
\sum_{J \succeq K \succeq I} 
\frac{(-1)^{ \ell(K)} \pi(I)}{\pi_u (J, K)\, \ell(K, I)} \rp \,  \Psi^J\, . 
\end{align*}
Hence, we get Eq.\,(\ref{Psi-Phi}). By a similar argument, 
(\ref{Phi-Psi}) follows by combining
Eq.\,(\ref{S-Psi}) with Eq.\,(\ref{Phi-S}).
\epfv

{\it
\begin{enumerate}
\item[$\bullet$] The relations between $\Lambda$ and $\Xi$: 
\begin{align}
\Lambda^I & = (-1)^{|I|} \sum_{J\succeq I}  \frac{(-1)^{\ell (J)}}{\pi_u( J, I)}\,\,  \Xi^J\, ,
\label{Xi-Lambda} \\
\Xi^I & = (-1)^{|I|} \sum_{J \succeq I} (-1)^{\ell(J)} lp (J, I) \,\,  \Lambda^J\, .
\label{Lambda-Xi}
\end{align}
\end{enumerate} }

\pf 
By applying the anti-involution $\omega_\Lambda$ 
in Proposition \ref{omega-Lambda} to Eq.\,$(\ref{Lambda-Psi})$,
and then, by Eq.\,$(\ref{omega-Lambda-e3})$ in the same proposition,   
we get 
\begin{align*}
\Xi^{\bar I} =(-1)^{\ell(I)} \sum_{ J \succcurlyeq  I}(-1)^{\ell(J)}
{fp (J, I)} \, \, \Lambda^{\bar J}\, . 
\end{align*}
Note that, for any compositions $I$ and $J$, 
$J \succcurlyeq  I$ iff $\bar J \succcurlyeq \bar I$. 
By replacing $\bar I$ by $I$ and $\bar J$ by $J$ in 
the equation above, 
we get

\begin{align*}
\Xi^{I}=(-1)^{\ell(\bar I)}
\sum_{ J \succcurlyeq  I}(-1)^{\ell(\bar J)}
{fp (\bar J, \bar I)}\, \, \Lambda ^{J}\, .
\end{align*}

Note that, $\ell(K)=\ell(\bar K)$ for any composition $K$.  
For any composition $I$ and $J$ with 
$J \succcurlyeq  I$, by Eqs.\,(\ref{lp-JI}) and (\ref{fp-JI}), 
it is easy to see that $lp (\bar J, \bar I)=fp(J, I)$.
With these observations and the equation above, we have
\begin{align*}
\Xi^{I}= (-1)^{\ell(I)}
\sum_{ J \succcurlyeq  I}(-1)^{\ell(J)}
{lp (J, I)}\, \, \Lambda ^{J}\, . 
\end{align*}

Hence, we get Eq.\,(\ref{Lambda-Xi}). 
Eq.\,(\ref{Xi-Lambda}) can be proved similarly
by applying $\omega_\Lambda$ to Eq.\,(\ref{Psi-Lambda}).
\epfv

The next two pairs of identities 
can also be proved similarly as above.

{\it
\begin{enumerate}

\item[$\bullet$] The relations between $S$ and $\Xi$:
\begin{align}
S^I & =  \sum_{J\succeq I}  \frac{1}{\pi_u( \bar{J}, \bar{I} )}\, \,  \Xi^J\, ,
\label{Xi-S} \\
\Xi^I & = (-1)^{\ell(I) }\sum_{J \succeq I} (-1)^{ \ell(J)} fp (J, I) \, \,  S^J\, . 
\label{S-Xi}
\end{align}

\item[$\bullet$] The relations between $\Xi$ and $\Phi$: 
\begin{align}
\Phi^I & =(-1)^{\ell(I) }
\sum_{J \succeq I} \lp 
\sum_{J \succeq K \succeq I} 
\frac{(-1)^{ \ell(K)} \pi(I)}{\pi_u (\bar{J}, \bar{K})\, \ell(K, I)} \rp \, \,  \Xi^J\, , 
 \label{Xi-Phi} \\
\Xi^I & = (-1)^{\ell(I) } \sum_{J \succeq I}
\lp \sum_{J \succeq K \succeq I}   \frac {(-1)^{\ell(K)}\, fp (K, I)}{sp (J, K)} \rp \, \, \Phi^J\, .
\label{Phi-Xi}
\end{align}
\end{enumerate} }

Finally, let us consider the relations between $\Psi$ and $\Xi$.

\begin{lemma}\label{L4.1.5}
For any composition $I$, we have 
\begin{align} 
\Xi^I &=\sum_{K\succcurlyeq I} c_{I, K} \,\,  \Psi^K\, ,  \label{Psi-Xi}\\
\Psi^I &=\sum_{K\succcurlyeq I} c_{\bar I, \bar K} \,\,  \Xi^K\, , \label{Xi-Psi}
\end{align}
where, for any composition $K \succcurlyeq I$, 
\begin{align} \label{c_IK}
c_{I, K}:=\sum_{K \succcurlyeq J \succcurlyeq I}
(-1)^{\ell(J)-\ell(I)} \, \, 
\frac {{fp (J, I)}}{\pi_u( K, J)} .
\end{align}
\end{lemma}

\pf Combining Eq.\,$(\ref{S-Xi})$ with Eq.\,(\ref{Psi-S}), 
we get
\begin{align*}
\Xi^{I} &= \sum_{ J \succcurlyeq  I}(-1)^{\ell(J)-\ell(I)}
{fp (J, I)} \sum_{K \succcurlyeq J}
\frac 1{\pi_u( K, J)} \, \, \Psi^K \\
&= \sum_{ K \succcurlyeq  I}
\left [ \sum_{K \succcurlyeq J \succcurlyeq I}
(-1)^{\ell(J)-\ell(I)} 
\frac {{fp (J, I)}}{\pi_u( K, J)} \right ] \, \, \Psi^K\, .
\end{align*}

Hence we get Eq.\,$(\ref{Psi-Xi})$. 
Eq.\,$(\ref{Xi-Psi})$ can be easily proved by 
applying the anti-involution $\omega_\Lambda$ 
in Proposition \ref{omega-Lambda}
to Eq.\,$(\ref{Psi-Xi})$ and then applying  
Eq.\,$(\ref{omega-Lambda-e3})$.
\epfv

Note that, we can also apply the involution 
$\tau$ in Proposition \ref{tau-flip}, 
instead of the anti-involution $\omega_\Lambda$, 
to Eq.\,$(\ref{Psi-Xi})$ to get another formula 
for $\Psi^I$ $(I\in \cC)$ in terms of $\Xi^J$ 
$(J \in \cC)$.

\begin{corol}
For any $I\in \cC$, we have
\begin{align} 
\Psi^I &=\sum_{K\succcurlyeq I} (-1)^{\ell(I)-\ell(K)} c_{I, K} \,\,  \Xi^K\, . 
\label{Xi-Psi-2}
\end{align}
\end{corol}

Furthermore, by comparing Eqs.\,$(\ref{Xi-Psi})$, 
$(\ref{Xi-Psi-2})$ and noting that the monomials $\Xi^K$ $(K\in \cC)$
are free in the $K$-algebra $\cNsf$, we get the following identity
for the coefficient $c_{I, K}$ for any  $I, K\in \cC$ 
with $K\succeq I$.
\begin{align} 
c_{\bar{I}, \bar{K}}= (-1)^{\ell(I)-\ell(K)} c_{I, K}.
\end{align}

More explicitly, combining with Eq.\,(\ref{c_IK}), 
it is easy to check that, 
for any $I, K\in \cC$ with $K\succeq I$, 
we have
\begin{align*} 
(-1)^{\ell(K)}
\sum_{K \succcurlyeq J \succcurlyeq I}
\frac {(-1)^{\ell(J)} \, fp (J, I)}{\pi_u( K, J)} 
=(-1)^{\ell(I)} 
\sum_{ \bar{K} \succcurlyeq J \succcurlyeq \bar{I} } 
\frac {(-1)^{\ell(J)} \,fp (J, \bar{I} )}{\pi_u( \bar{K}, J)}.
\end{align*}

\subsection{Formulas for $e^{-u\Phi(t)}$} \label{S4.2}

Let $u$ be another central parameter, i.e. 
it commutes with $t$ and any NCSF's in $\cNsf$. 
In this section, we derive some formulas 
for $e^{-u\Phi(t)}$ in terms of the NCSF's 
in the universal \cNcs system $(\cNsf, \Pi)$. 
These formulas later will be needed in 
Subsection \ref{S5.2} 
for the study of the formal
flows generated by 
$F_t\in \ataz$.

Let us start with the following two lemmas.
 
\begin{lemma}\label{Flow-lemma}
\begin{align}
e^{-u\Phi(t)}&=1+\sum_{I\in \cC } \frac{(-u)^{\ell(I)} t^{|I|}}{sp(I)} \Phi^I\, ,
\label{Flow-lemma-e0}   \\
e^{-u\Phi(t)}&=1+ \sum_{I \in \mathcal C} (-1)^{|I|-\ell(I)} t^{|I|} 
\lp \sum_{I\succeq J}
\frac{(-1)^{\ell(J)} u^{ \ell(J)} }{\ell(J)!\,  \ell (I, J)} \rp 
\,  \Lambda^I  \, , \label{Flow-lemma-e1} \\
e^{-u\Phi(t)}&= 1+ \sum_{I \in \mathcal C} (-1)^{\ell(I)} t^{|I|} 
\lp \sum_{I\succeq J}
\frac{u^{\ell(J)}}{\ell(J)!\,  \ell (I, J)} \rp 
\,  S^I  \, ,\label{Flow-lemma-e2}  \\
e^{-u\Phi(t)}&= 1+ \sum_{I \in \mathcal C}  t^{|I|}
\lp 
\sum_{I \succeq K \succeq J} 
\frac{(-1)^{ \ell(K)} u^{\ell(J)} }{\pi_u (I, K)\, \ell(K, J) \ell(J)!} \rp \,  
\Psi^I  \, ,\label{Flow-lemma-e3}  \\
e^{-u\Phi(t)}&= 1+ \sum_{I \in \mathcal C}  t^{|I|}
\lp 
\sum_{I \succeq K \succeq J} 
\frac{(-1)^{ \ell(K)} u^{\ell(J)}}{\pi_u (\bar{I}, \bar{K})
\, \ell(K, J)\, \ell(J)!} \rp \,  \Xi^I \, . \label{Flow-lemma-e4}
\end{align}
\end{lemma}

\pf First, by Eq.\,(\ref{Phi(t)}), we have
\allowdisplaybreaks{
\begin{align}
e^{-u\Phi(t)}&=1+\sum_{m\geq 1} \frac {(-u)^m}{m!}
(\sum_{k\geq 1} t^k \frac{\Phi_k}k )^m \\
&=1+\sum_{m\geq 1} \frac {(-u)^m}{m!} 
\sum_{\substack{I\in \cC\\ \ell(I)=m} } \frac {t^{|I|}} {\pi(I)}\, \Phi^I \nno \\
&=1+\sum_{I\in \cC } \frac{(-u)^{\ell(I)} t^{|I|}}{\ell(I)\pi(I)}\, \Phi^I \nno \\
&=1+\sum_{I\in \cC } \frac{(-u)^{\ell(I)} t^{|I|}}{sp(I)}\, \Phi^I \, .\nno
\end{align} }

Therefore we get Eq.\,(\ref{Flow-lemma-e0}).
All other formulas in the lemma follow from 
Eq.\,(\ref{Flow-lemma-e0}) 
and the relations of involved NCSF's with $\Phi$. 
As one example, we give a proof for Eq.\,(\ref{Flow-lemma-e1}). 
The proofs for Eqs.\,(\ref{Flow-lemma-e2})--(\ref{Flow-lemma-e4}) are similar. 

First, by Eqs.\,(\ref{Flow-lemma-e0}) and (\ref{Lambda-Phi}), 
we have
\begin{align*}
e^{-s\Phi(t)}&= 1+ \sum_{I\in \mathcal C} 
\frac{ (-1)^{|I|-\ell(I)} t^{|I|} u^{\ell(I)} }{sp(I)}\,  
\sum_{J\succeq I} 
\frac{ (-1)^{ \ell(J)} \pi(I) }{\ell (  J, I)}\, \,  \Lambda^J  \\
\intertext{Switching the order of the summations and noting that $|I|=|J|$:}
&=1 + \sum_{J \in \mathcal C} (-1)^{|J|-\ell(J)} t^{|J|} 
\lp \sum_{J\succeq I}
\frac{(-1)^{ \ell(I)} u^{\ell(I)}}{\ell(I)!\,  \ell (J, I)} \rp 
\,  \Lambda^J  \\
\intertext{Switching the summation indices $I$ and $J$:}
&=1+ \sum_{I \in \mathcal C} (-1)^{|I|-\ell(I)} t^{|I|} 
\lp \sum_{I\succeq J}
\frac{(-1)^{ \ell(J)} u^{\ell(J)}}{\ell(J)!\,  \ell (I, J)} \rp 
\,  \Lambda^I \, .
\end{align*}
\epfv

\begin{lemma}\label{Flow-lemma2}
\begin{align}
e^{-\Phi(t)} & = 1 + \sum_{m\geq 1}(-1)^m t^m \, \, \Lambda_m \, , 
\label{Flow1-Lambda} \\
e^{ \Phi(t)} & = 1 + \sum_{m\geq 1} t^m \sum_{I\in \cC_m} 
\frac 1{\pi_u(I)} \, \, \Psi^I , \label{Flow11-Psi} 
\end{align}
\end{lemma}
\pf Note that, by Eqs.\,(\ref{UE-1}), (\ref{UE-2}) for 
the universal \cNcs system $(\cNsf, \Pi)$ and 
Eqs.\,(\ref{lambda(t)}) and (\ref{sigma(t)}), 
we have 
\begin{align}
e^{-\Phi(t)}&=\sigma(t)^{-1}=\lambda(-t)=1+\sum_{m\geq 1} (-1)^m t^m \Lambda_m\, ,
\label{Flow-lemma2-pe1} \\
 e^{\Phi(t)}&=\sigma(t)=\lambda(-t)=1+\sum_{m\geq 1}  t^m S_m\, .
 \label{Flow-lemma2-pe2}
\end{align}

Therefore, we have Eq.\,$(\ref{Flow1-Lambda})$. 
Furthermore, Eq.\,$(\ref{Psi-S})$ 
with $I=\{m\}$, we have
\begin{align}
S_m=\sum_{I\in \cC_m} 
\frac 1{\pi_u(I)}. 
\end{align}
Combining the equation above with 
Eq.\,$(\ref{Flow-lemma2-pe2})$, we get 
Eq.\,$(\ref{Flow11-Psi})$.
\epfv


Now we can formulate the main result of this subsection as follows.

\begin{propo}\label{Flow-propo}
\begin{align}
e^{-u\Phi(t)} &=1+ \sum_{I \in \mathcal C} (-1)^{|I|} t^{|I|} 
\binom{u}{\ell(I)}\,  \Lambda^I  \, , \label{Flow-propo-e1} \\
e^{-u\Phi(t)} &= 1+ \sum_{I \in \mathcal C}  t^{|I|} 
\binom{-u}{\ell(I)}\, 
\,  S^I \, ,\label{Flow-propo-e2}  \\
e^{-u\Phi(t)} &= 1+ \sum_{I \in \mathcal C}  t^{|I|}
\lp 
\sum_{I \succeq J} 
\frac{ \binom{-u}{\ell(J)} }{\pi_u (I, J)} \rp \,  
\Psi^I  \, ,\label{Flow-propo-e3}  \\
e^{-u\Phi(t)} &= 1+ \sum_{I \in \mathcal C} (-1)^{\ell(I)} t^{|I|}
\lp 
\sum_{I \succeq J} 
\frac{ \binom{u}{\ell(J)} }{\pi_u (\bar{I}, \bar{J} )} \rp \, 
  \Xi^I  \, . \label{Flow-propo-e4}
\end{align}
\end{propo}

\pf Let us first show Eq.\,(\ref{Flow-propo-e1}). 
From Eq.\,(\ref{Flow-lemma-e1}), 
we see that there exist (unique) polynomials 
$\varphi_{{}_\Lambda} (I, u)\in K[u]$ $(I\in \cC)$ of 
degree $d\leq \ell(I)$ 
such that
\begin{align}
e^{-u \Phi(t)} = 1+ & \sum_{I\in \cC} t^{|I|} \varphi_{{}_\Lambda} (I, u)\, \Lambda^I\, ,
\label{Flow-propo-pe2}\\
\varphi_{{}_\Lambda} (I, 0)&=0 \, ,\label{u-0}
\end{align}
for any $I\in \cC$.

By the fact $e^{-(u+1) \Phi(t)}=e^{- \Phi(t)}e^{-u \Phi(t)}$ 
and Eq.\,(\ref{Flow-propo-pe2}) above, we have
\begin{align*}
& 1+\sum_{I\in \cC}t^{|I|} \varphi_{{}_\Lambda} (I, u+1) \Lambda^I 
 \\ & \qquad \quad
=( 1+\sum_{I\in \cC} t^{|I|} \varphi_{{}_\Lambda} (I, 1) \Lambda^I)
( 1+\sum_{I\in \cC} t^{|I|} \varphi_{{}_\Lambda} (I, u) \Lambda^I ) \, .
\end{align*} 
By comparing the coefficients of $\Lambda^I$ $(I\in \cC)$ and 
noting that they are free in the $K$-algebra $\cNsf$, 
it is easy to see that, 
for any $I\in \cC$, we have
\begin{align}\label{Delta-Relation}
 \Delta ( \varphi_{{}_\Lambda} (I, u) )
=  \varphi_{{}_\Lambda} (I, 1)+
\sum_{\substack{J, K\in \cC \\ J\cdot K=I}}  
\varphi_{{}_\Lambda} (J, 1) \, \varphi_{{}_\Lambda} (K, u),
\end{align}
where $\Delta: K[u]\to K[u]$ is the {\it difference operator} 
which maps any $q(u)\in K[u]$ to $q(u+1)-q(u)$. 

Recall that, we have the following well-known facts. First, 
for any $m\geq 0$, we have 
\begin{align}\label{Flow-propo-pe4}
\Delta \binom{u}{m}=\binom{u}{m-1}.
\end{align}
Secondly, for any polynomial $q(u)\in K[u]$ 
of degree $d\geq 0$, we have
\begin{align}\label{Flow-propo-pe5}
q(u)=\sum_{m = 0}^d  a_m \binom{u}{m},
\end{align}
where $a_m=(\Delta^m q) (0)$ for any $0\leq m\leq d$.

Now, we apply the facts above to the polynomials 
$\varphi_{{}_\Lambda}(I, u)\in K[u]$ $(I\in \cC)$. 
By comparing Eq.\,(\ref{Flow1-Lambda}) with 
Eq.\,(\ref{Flow-propo-pe2}) with $u=1$, and again noting that
the monomials $\Lambda^I$ $(I\in \cC)$ are free in $\cNsf$, 
we have
\begin{align*}
\varphi_{{}_\Lambda} (I, 1)=
\begin{cases}
(-1)^{|I|} & \text{\quad if \quad} \ell(I)=1, \\
0& \text {\quad if \quad} \ell(I)\geq 2.
\end{cases}
\end{align*}

By the equation above and Eqs.\,(\ref{u-0}), (\ref{Delta-Relation}), 
it is easy to check that, for any $I\in \cC$ and $m\geq 0$, 
we have
\begin{align*}
\left. \Delta^m \varphi_{{}_\Lambda}  (I, u)\, \right |_{u=0} 
=
\begin{cases}
0  & \text{\quad if \quad} m \neq \ell(I), \\
(-1)^{|I|} & \text{\quad if \quad} m=\ell(I). 
\end{cases}
\end{align*}
Then, by the general fact given by
Eq.\,(\ref{Flow-propo-pe5}), we have, for any $I\in \cC$,
\begin{align*}
\varphi_{{}_\Lambda} (I, u)=(-1)^{|I|} \binom{u}{\ell(I)}.
\end{align*}
Combining the equation above with Eq.\,(\ref{Flow-propo-pe2}), 
we get Eq.\,(\ref{Flow-propo-e1}).

To show Eq.\,(\ref{Flow-propo-e2}), we first apply the involution 
$\tau:\cNsf\to \cNsf$ in Proposition \ref{tau-flip} to Eq.\,(\ref{Flow-propo-e1}).
By Eq.\,(\ref{tau-flip-e1}), we have
\begin{align*}
e^{u\Phi(t)} &= 1+ \sum_{I \in \mathcal C}  
t^{|I|} \binom{u}{\ell(I)}\, 
\,  S^I \, 
\end{align*}
Then, replacing $u$ by $-u$ in the equation above, 
we get  Eq.\,(\ref{Flow-propo-e2}).

Next, we show Eq.\,(\ref{Flow-propo-e3}). 
Note that, by applying the involution $\tau$ 
in Proposition \ref{tau-flip} to 
Eq.\,(\ref{Flow-propo-e3}), we will get 
Eq.\,(\ref{Flow-propo-e4}).

First, by Eq.\,(\ref{Flow-lemma-e3}) 
with $u$ replaced by $-u$, 
there exist (unique) polynomials 
$\varphi_{{}_\Psi} (I, u)\in K[u]$ $(I\in \cC)$ of 
degree $d\leq \ell(I)$ 
such that
\begin{align}
e^{u \Phi(t)} = 1+& \sum_{I\in \cC} t^{|I|} \varphi_{{}_\Psi} (I, u) \Lambda^I\, ,
\label{Flow-propo-pe7} \\
\varphi_{{}_\Psi} (I, 0)&=0 \, ,\label{3u-0}
\end{align}
for any $I\in \cC$.

By the fact $e^{(u+1) \Phi(t)}=e^{\Phi(t)}e^{u \Phi(t)}$ 
and Eq.\,(\ref{Flow-propo-pe7}) above, 
we see that the polynomials 
$\varphi_{{}_\Psi} (I, u)\in K[u]$ 
$(I\in \cC)$ also satisfy 
Eq.\,(\ref{Delta-Relation}).
On the other hand, by Eqs.\,(\ref{Flow11-Psi}) and 
(\ref{Flow-propo-pe7}) with $u=1$, 
and also the freeness of 
the monomials $\Psi^I$ $(I\in \cC)$, 
we have, for any $I\in \cC$,
\begin{align}
\varphi_{{}_\Psi} (I, 1)=
\frac{1}{\pi_u(I)}. \label{Flow-propo-pe8}
\end{align}

Then, by Eqs.\,(\ref{3u-0}), (\ref{Flow-propo-pe8}) and 
Eq.\,(\ref{Delta-Relation}) for $\varphi_{{}_\Psi} (I, u)$, 
it is easy to see that, for any $I\in \cC$ and $m\geq 0$, 
we have
\begin{align*}
\left. \Delta^m \varphi_{{}_\Psi}(I, u)\, \right |_{u=0} 
=\sum_{\substack {I\succeq J\\ \ell(J)=m }}
\frac{1}{\pi_u(I, J)}.
\end{align*}

By the general fact given by 
Eq.\,(\ref{Flow-propo-pe5}), we have, for any $I\in \cC$,
\begin{align*}
\varphi_{{}_\Psi} (I, u)= \sum_{I\succeq J}
\frac{\binom{u}{\ell(J)}}{\pi_u(I, J)}.
\end{align*}

Hence, we have
\begin{align*}
e^{u\Phi(t)} &= 1+ \sum_{I \in \mathcal C}  t^{|I|} 
\sum_{I\succeq J}
\frac{\binom{u}{\ell(J)}}{\pi_u(I, J)} \,  \Psi^I \, .
\end{align*}
Replacing $u$ by $-u$ in the equation above, we get 
Eq.\,(\ref{Flow-propo-e3}).
\epfv

Finally, by comparing the formulas Eqs.\,(\ref{Flow-lemma-e1}) 
and (\ref{Flow-propo-e1})
and using the freeness of the NCSF's $\Lambda^K$ $(K\in \cC)$, 
it is east to see that 
we have the following identity for composition.

\begin{corol}
For any composition $I\in\cC$ and a free variable $u$, 
we have
\begin{align}
 \sum_{I\succeq J}
\frac{(-1)^{\ell(J)} u^{ \ell(J)} }{\ell(J)!\,  \ell (I, J)} 
  =(-1)^{\ell(I)} \binom{u}{\ell(I)}\, . 
\end{align}
\end{corol}

\renewcommand{\theequation}{\thesection.\arabic{equation}}
\renewcommand{\therema}{\thesection.\arabic{rema}}
\setcounter{equation}{0}
\setcounter{rema}{0}

\section{\bf Applications to the Inversion Problem}\label{S5}

In this section, 
we mainly apply the identities of NCSF's derived 
in the previous section and
the specialization 
$\cS_{F_t}$ $(F_t\in \ataz)$ 
in Theorem \ref{S-Correspondence}
to derive some formulas
for the inverse maps, the D-Log's and 
the formal flow of the automorphism 
$F_t\in \ataz$. 
In Subsections \ref{S5.1} and \ref{S5.2}, 
we derive formulas for the D-Log of and 
the formal flow generated by $F_t\in \ataz$, 
respectively, in terms of the differential 
operators in the \cNcs system $\Oft$.
In Subsection \ref{S5.3},
we mainly apply the results in the previous 
two subsections to the special automorphisms 
$F_t(z)\in \etaz$ to 
derive some inversion formulas. 
In Subsection \ref{S5.4}, motivated 
by the fact that 
$C_m(z) \in \kzz^{\times n}$
$(m\geq 1)$ in the commutative case
capture the nilpotence 
of the Jacobian matrix $JH$, 
we give formulas for 
$C_m(z)$'s in terms of 
the differential operators in $\Oft$ 
other than $\psi_m$'s.
Finally, in Subsection \ref{S5.5}, 
we discuss a connection of NCSF's with 
the well-known Jacobian conjecture.

\subsection{D-Log's in Terms of 
Other Differential Operators in $\Oft$}\label{S5.1}

Considering the important role played 
by the D-Log's in the inversion 
problem (see \cite{E1}--\cite{E3}, 
\cite{Z-exp} and \cite{WZ} 
for more discussions in the commutative case), 
we consider the expressions 
of the {\it D-Log} $a_t(z)$ of $F_t(z)$ 
(see page \pageref{D-Log})
in terms of other differential 
operators in the \cNcs  system $\Oft$.
Note that, the problem to express $\phi$ 
in terms of 
the derivations $\psi$ can be viewed 
as a special case of the so-called 
{\it the problem of continuous Baker-Campbell-Hausdorff 
exponents} in the mathematical physics 
(see $\S 4.10$ of \cite{G-T} 
and the references given there).

Recall that, by Eq.\,$(\ref{Def-d(t)})$ and 
the relation of $d(t)$ with the D-Log 
(see page \pageref{D-Log}), we have 
\begin{align}\label{D-Log-review}
d(t):=-\lb a_t(z)\pz \rb=\sum_{m=1}^\infty \frac {t^{m}}m \,\, \phi_m \, .
\end{align}

\begin{propo}\label{D-Log-Others}
\begin{align}
a_t(z) & =   \sum_{I\in \cC}  
 \frac{ (-1)^{ \ell(I)-|I|+1}  }{\ell (I)}\, \, t^{|I|}\, \lambda^I z\, , \label{D-Log-lambda}          \\
a_t(z) & =  \sum_{I\in \cC}  
 \frac{(-1)^{ \ell(I)}}{\ell (I)} \, \,  t^{|I|}\,  s^I  z\, ,           \\
a_t(z) & =  
\sum_{I\in \cC} 
\lp \sum_{ I \succeq J }   
\frac{(-1)^{ \ell(I)}}{\pi_u(I, J) \, \ell (J)} \rp \, t^{|I|}\, \psi^I z\, , \label{D-Log-psi} \\
a_t(z) & =  
\sum_{I\in \cC} 
\lp \sum_{ I \succeq J }   
\frac { (-1)^{\ell(J) } }{\pi_u (\bar{I}, \bar{J} )\, \ell(J)} \rp \, \, t^{|I|}\, \xi^I  z\, . 
\label{D-Log-xi}
\end{align}
\end{propo}

\pf
All the formulas above follow directly from the 
identities of $\phi$ with the corresponding differential operators.
For example, by Eq.\,(\ref{Xi-Phi}) with $I=\{m\}$, 
we have, for any $m\geq 1$,
\begin{align*} 
\Phi_m & =(-1)
\sum_{J \in \cC_m} \lp 
\sum_{J \succeq K } 
\frac{(-1)^{ \ell(K)} m }{\pi_u (\bar{J}, \bar{K} )\, \ell(K)} \rp \,  \Xi^J\\ 
& =(-1)
\sum_{I\in \cC_m} \lp 
\sum_{I \succeq J } 
\frac{(-1)^{ \ell(J)} m }{\pi_u (\bar{I}, \bar{J} )\, \ell(J)} \rp \,  \Xi^I. 
\end{align*}
Then, applying the specialization $\cSft$ to the equation above, 
and by Eq.\,(\ref{D-Log-review}), 
\begin{align*}
\lb a_t(z)\pz \rb  = -d(t)= \sum_{I\in \cC} 
\lp \sum_{ I \succeq J }   
\frac { (-1)^{\ell(J) } }{\pi_u (\bar{I}, \bar{J} )\, \ell(J)} \rp \, \, t^{|I|}\, \xi^I \, . 
\end{align*}
Applying the equation above to $z$, 
we get Eq.\,(\ref{D-Log-xi}). 
\epfv

Applying the specialization $\cS_{F_t}$ to 
Eq.\,$(71)$ 
in \cite{G-T} and then 
applying Corollary \ref{S-Correspondence-2}, 
we get the following improved formula for $d(t)$.

\begin{theo} \label{T4.1.7}
\begin{align} \label{T4.1.7-e1}
d(t)&=\sum_{r\geq 1} \int_0^t d t_1 \int_0^{t_1} dt_2 \cdots 
\int_0^{t_{r-1}} dt_r \\
&\quad   \quad\quad
\sum_{\sigma\in {\bf S}_r} 
\frac{(-1)^{d(\sigma)}}{r} \binom{r-1}{d(\sigma)}^{-1} 
h(t_{\sigma(r)}) \cdots h(t_{\sigma(1)})\, ,\nno
\end{align}
where ${\bf S_r}$ is the symmetric group of degree $r$ 
and, for any  $\sigma\in {\bf S}_r$, 
$d(\sigma)$ is the numbers of the descents of $\sigma$.
\end{theo}

Note that, by applying the equation above to $-z$, 
we get another formula for the D-Log $a_t(z)$ of $F_t$.

\subsection{Formal Flows in Terms of the Differential 
Operators in $\Oft$} \label{S5.2}

In this subsection, we consider the expressions of 
the formal flows, which has been studied in 
\cite{E1}--\cite{E3} and \cite{WZ} 
for the commutative case, in terms 
of the differential operators 
in the \cNcs system $\Oft$.

Let $u$ be another central parameter, i.e. 
it commutes with $z$ and $t$. We define
\begin{align}
F_t(z, u):= e^{u \lb a_t(z)\pz \rb }\,z= e^{-u \, d(t) }\,z. 
\label{Def-flow}
\end{align}

Note that, since $d(0)=0$, the exponential above 
is always well-defined. Actually it is easy to see 
$F_t(z, u)\in (K[u][[t]])
\langle \langle z \rangle \rangle^{\times n}$.
Therefore, for any $u_0 \in K$, $F_t(z, u_0)$ 
makes sense. 

Following its analog in 
\cite{E1}--\cite{E3} and \cite{WZ} 
in the commutative case, 
we call $F_t(z, u)$ the {\it formal flow} 
generated by $F_t(z)$ or simply 
the {\it formal flow} of $F_t(z)$.

Two remarks on the formal flows defined above 
are as follows.

First, it is well-known that 
the exponential of a derivation
of any $K$-algebra $A$, when it makes sense,   
is always 
an automorphism of the algebra, so in our case, 
for any $u_0\in K$, $e^{u_0 \lb a_t(z)\pz \rb }$ 
is also an automorphism of $\kttzz$ over $K[[t]]$ 
which maps $z$ to $F_t(z, u_0)$. 
From Eq.\,$(\ref{Def-flow})$, 
it is clear that this automorphism also lies
in $\ataz$ since $o(a_t(z)) \geq \alpha$.

Secondly, by Eq.\,$(\ref{Def-flow})$ 
and the remark above, 
the formal flow $F_t(z, u)$ has 
the following properties:
\begin{align}
F_t(z, 0)& = z, \\
F_t(z, 1)&=F_t(z),\\
F_t( F_t(z, u_2), u_1)&=F_t(z, u_1+u_2),
\end{align}
for any $u_1, u_2\in K$.

In other words, $F_t(z, u)$ forms an one-parameter 
subgroup of the group $\ataz$. 
Therefore, for any integer $m\in \bZ$, 
$F_t(z, m)$ gives the $m^{th}$ 
(composing)\label{composingPower} power of $F_t$ 
as an element of the group $\ataz$. 
In particular, by setting $m=-1$, 
we get the inverse map $G_t$ of $F_t$, 
i.e. $F_t(z, -1)=G_t(z)$.

The main result of this subsection is the following proposition 
which expresses the D-Log of $F_t\in \ataz$ in terms of 
the differential operators in the \cNcs $\Oft$. 

\begin{propo}\label{flow-propo}
\begin{align}
F_t(z, u)&=z+ \sum_{I \in \mathcal C} \frac {(-u)^{\ell(I)}}{sp(I)}
\,  t^{|I|}\,\phi^I z \, , \label{flow-propo-e0} \\
F_t(z, u)&=z+ \sum_{I \in \mathcal C} (-1)^{|I|}
\binom{u}{\ell(I)}\,  t^{|I|}\, \lambda^I z \, , \label{flow-propo-e1} \\
F_t(z, u)&= z+ \sum_{I \in \mathcal C}  
\binom{-u}{\ell(I)}\, 
\, t^{|I|}\, s^I z \, ,\label{flow-propo-e2}  \\
F_t(z, u)&= z+ \sum_{I \in \mathcal C}  
\lp 
\sum_{I \succeq J} 
\frac{ \binom{-u}{\ell(J)} }{\pi_u (I, J)} \rp \,  
t^{|I|}\, \psi^I  z\, ,\label{flow-propo-e3}  \\
F_t(z, u)&= z+ \sum_{I \in \mathcal C} (-1)^{\ell(I)} 
\lp 
\sum_{I \succeq J} 
\frac{ \binom{u}{\ell(J)} }{\pi_u (\bar{I}, \bar{J} )} \rp \, 
  t^{|I|}\, \xi^I  z\, . \label{flow-propo-e4}
\end{align}
\end{propo}

\pf 
By the expression of $d(t)$ given in 
Eq.\,(\ref{D-Log-review}) and the correspondence 
in Corollary \ref{S-Correspondence-2}, 
we see that the formal flow is given by
\begin{align}
F_t(z, u)= \cSft(e^{-u\Phi(t)})\, z.
\end{align}
Then, by the equation above, it is easy to see that,
Eq.\,(\ref{flow-propo-e0}) follows directly from 
Eq.\,(\ref{Flow-lemma-e0}) and 
Eqs.\,$(\ref{flow-propo-e1})$--$(\ref{flow-propo-e4})$
follows respectively from 
Eqs.\,$(\ref{Flow-propo-e1})$--$(\ref{Flow-propo-e4})$.
\epfv

\subsection{Some Inversion Formulas}\label{S5.3}
In this subsection, we mainly apply the identities in 
Subsection \ref{S4.1} to derive some inversion formulas. 

First, let us consider 
the Taylor series expansions of 
$u(F_t)$ and $u(G_t)$ for any 
$F_t\in \ataz$ and $u(z)\in \kzz$ 
in terms of the differential operators 
$\lambda$ and $\psi$.

\begin{propo}\label{Taylor-lambda-psi}
For any $F_t\in \ataz$ $(\alpha\geq 1)$ and $u(z)\in \kzz$, we have
\begin{align}
u(F_t(z))& = u(z)+ \sum_{m\geq 1}(-1)^m t^m \, \, \lambda_m \, u(z), 
\label{Taylor-lambda-psi-e1} \\
u(F_t(z))& = u(z) + \sum_{m\geq 1} t^m \sum_{I\in \cC_m} 
\frac {(-1)^{\ell (I)}} { \pi_u(\bar{I}) }  \, \, \psi^I  u(z)
 \label{Taylor-lambda-psi-e2} \\
u(G_t(z))& = u(z)+ \sum_{m\geq 1} t^m \sum_{I\in \cC_m} 
(-1)^{\ell(I)-m} \, \, \lambda^I  u(z), \label{Taylor-lambda-psi-e3} \\
u(G_t(z)) & = u(z)+ \sum_{m\geq 1} t^m \sum_{I\in \cC_m} 
\frac 1{\pi_u(I)} \, \, \psi^I  u(z), \label{Taylor-lambda-psi-e4}
\end{align}
\end{propo}
\pf 
Eq.\,(\ref{Taylor-lambda-psi-e1}) follows from the definition 
of the differential operators $\lambda_m$'s 
(see Eqs.\,$(\ref{NewTaylorExpansion-e1})$ 
and $(\ref{Def-f(t)})$). 
Eq.\,(\ref{Taylor-lambda-psi-e2}) follows from 
Eqs.\,$(\ref{Taylor-lambda-psi-e1})$ and 
$(\ref{Psi-Lambda})$.
Eq.\,(\ref{Taylor-lambda-psi-e3}) follows from 
Eqs.\,$(\ref{NewTaylorExpansion-e2})$,  $(\ref{Def-g(t)}$) 
and $(\ref{Lambda-S})$. Finally, 
Eq.\,(\ref{Taylor-lambda-psi-e4}) follows from 
Eqs.\,$(\ref{NewTaylorExpansion-e2}$),  
$(\ref{Def-g(t)})$ and $(\ref{Psi-S})$.
\epfv

{\it From now on and throughout 
the rest of this paper, 
we will assume $F_t(z)\in \etaz$, i.e. 
$F_t(z)$ is an automorphism 
of the form $F_t(z)=z-tH(z)$ for some 
$H(z)\in \kzz^{\times n}$ with 
$o(H(z))\geq \alpha$. 
We will also freely use 
the notations fixed 
in Section \ref{S3}.}

\vskip2mm

First, by Eq.\,$(\ref{Taylor-lambda-psi-e3})$ with $u(z)=z$, 
we get the following inversion formula for $F_t(z)$. 

\begin{propo}\label{Inv-lambda}
For any $m\geq 1$, we have
\begin{align}\label{Inv-lambda-e1}
N_{[m]}(z)&=(-1)^{m}\sum_{I\in \cC_m}(-1)^{\ell(I)} \, \, \lambda^I   z ,
\end{align}
where $\lambda=\{ \lambda_k \,|\,k \geq 1\}$ are given by 
Eq.\,$(\ref{Special-lambda-2})$ in general and by Eq.\,$(\ref{Special-lambda-3})$ 
when $z$ are commutative free variables.
\end{propo}

Applying Eq.\,$(\ref{Taylor-lambda-psi-e4})$ with $u(z)=z$, we get 
the following inversion formula in terms of $\psi$. 

\begin{propo}\label{Inv-psi}
For any $m\geq 1$, we have
\begin{align}\label{Inv-psi-e1}
N_{[m]}(z)=\sum_{I\in \cC_m} \frac 1{\pi_u(I)} \, \, \psi^I  z\, ,
\end{align}
where $\psi=\{ \psi_k \,|\,k \geq 1\}$ are given by 
Eq.\,$(\ref{Special-psi-1})$ in general and by Eq.\,$(\ref{Special-psi-m})$ 
when $z$ are commutative free variables.
\end{propo}

By using the identities between the NCSF's $\Xi$ 
and $\Psi$, we can get another inversion formula 
in terms of $\psi$ as follows.

\begin{propo} \label{Inv-psi-II}
For any $m\geq 1$, we have 
\begin{align}\label{Inv-psi-2}
N_{[m]}(z)=\sum_{I\in \cC_m} c_I \, \psi^I  z,
\end{align}
where, for any composition $I$, $c_I$ is given by 
\begin{align}\label{c_I}
c_I =\sum_{I \succcurlyeq J } (-1)^{\ell (J)-1} \frac {fp(J)}{\pi_u (I, J)}.
\end{align}
\end{propo}

\pf First, for any fixed $m\geq 1$, let $I$ be the composition 
$\{m \}$ (of length $1$). For any composition $J$, we have, 
$J\succcurlyeq I$ iff $|J|=m$, and in this case, 
$fp(J, I)=fp(J)$ by Eq.\,(\ref{fp-JI}). 
With these observations and Eq.\,(\ref{c_IK}), 
it is easy to see that $c_{I, K}=c_K$
for any $K\in \cC_m$ and 
Eq.\,$(\ref{Psi-Xi})$
becomes 
\begin{align*}
\Xi_m=\sum_{K\in \cC_m} c_K \, \Psi^K= 
\sum_{I\in \cC_m} c_I \, \Psi^I\, .
\end{align*}

Applying $\cS_{F_t}$ to the equation above and 
then applying the resulted equation to $z$, by Eq.\,(\ref{Special-xi}),
we get Eq.\,$(\ref{Inv-psi-2})$.
\epfv

Next, let us derive the following recurrent inversion formula.

\begin{propo}\label{Recur-Inv}
We have the following recurrent inversion formula.
\begin{align}
N_{[1]}(z) &=H (z), \label{Recur-Inv-e1} \\
N_{[m]}(z) &=\frac{m}{m-1} \sum_{\substack{I\in \cC_m \\ \ell (I)\geq 2}} 
\frac 1{\pi_u(\bar I)}
\lb N_{[i_1]}\pz\rb  
\cdots \lb N_{[i_{k-1}]} \pz \rb  N_{[i_k]}(z), \label{Recur-Inv-e2}
\end{align}
for any $m\geq 2$.
\end{propo}

For a different but more effective recurrent inversion formula, 
see Theorem $5.5$ in \cite{NC-IVP}.

\pf First, Eq.\,$(\ref{Recur-Inv-e1})$ 
follows easily from 
Eqs.\,(\ref{Inv-psi-e1}), (\ref{Special-psi-1}) with $m=1$
and Eq.\,$(\ref{L4.1.1-e3})$. 
 
To show Eq.\,$(\ref{Recur-Inv-e2})$,  
for any $m\geq 1$,  by Eq.\,$(\ref{Xi-S})$ with $I=\{ m \}$, 
we have
\begin{align} \label{Recur-Inv-pe1}
S_m =\sum_{I\in \cC_m} \frac 1{\pi_u( \bar I )}\,\, \Xi^{I}\, . 
\end{align}
Then, we apply the specialization 
$\cS_{F_t}$ to Eq.\,$(\ref{Recur-Inv-pe1})$ above and, by  
Corollary \ref{S-Correspondence-2}, 
we get
\begin{align}\label{Recur-Inv-pe2}
s_m =\sum_{I\in \cC_m} \frac 1{\pi_u(\bar I)}\, \xi^I.
\end{align}

Note that, by Eq.\,$(\ref{NewTaylorExpansion-e2})$ with $u(z)=z$ 
and Eqs.\,$(\ref{Def-g(t)})$, $(\ref{Special-xi})$, we have
$s_m   z=\xi_m\cdot z=N_{[m]}(z)$ for any $m\geq 1$.
Then, by applying  both sides of
Eq.\,$(\ref{Recur-Inv-pe2})$ to $z$, we get
\begin{align*}
N_{[m]}(z) &=\sum_{I\in \cC_m} \frac 1{\pi_u(\bar I)}\, \xi^I \,  z, \\
N_{[m]}(z) &=\frac 1m N_{[m]}(z) + \sum_{\substack{I\in \cC_m \\ \ell(I)\geq 2}} 
\frac 1{\pi_u(\bar I)}\, \xi^I  \,  z,\\
\frac {m-1}m & N_{[m]}(z)= \sum_{\substack{I\in \cC_m \\ \ell(I)\geq 2}} 
\frac 1{\pi_u(\bar I)}\, \xi^I  \,  z.
\end{align*}
Therefore, we have
\begin{align*}
N_{[m]}(z) &= \frac{m}{m-1} \sum_{\substack{I\in \cC_m \\ \ell(I)\geq 2}} 
\frac 1{\pi_u(\bar I)}\, \xi^I  \,  z \\
\intertext{Applying Eq.\,$(\ref{Special-xi})$:}
&=\frac{m}{m-1} \sum_{\substack{I\in \cC_m \\ \ell (I)\geq 2}}  \frac 1{\pi_u( \bar I)}
\lb N_{[i_1]}\frac{\p}{\p z}\rb 
\cdots \lb N_{[i_{k-1}]}\frac{\p}{\p z}\rb  N_{[i_k]}(z).
\end{align*}
\epfv

Finally, let us end this subsection with 
the following identity of differential operators, 
which does not seem to be obvious.

\begin{propo}\label{Identity-H(z)-pz}
For any $m\geq 1$ and $H(z) \in \kzz^{\times n}$, 
let $\psi_m$ be the derivation given 
by Eq.\,$(\ref{Special-psi-1})$. Then, we have

\begin{align}\label{Identity-H(z)-pz-e1}
\frac 1{m!} \left. \lb  H(w) \pz \rb^m \,  \right |_{w=z}
= \sum_{I\in \cC_m}
\frac {(-1)^{\ell (I)}} { \pi_u(\bar{I}) } \, \, \psi^I \, . 
\end{align}

In particular, when $z$ are commutative free variables, 
we have
\begin{align}
& \psi_m = \lp (JH)^{m-1}H \rp \pz\, , \label{Identity-H(z)-pz-e2} \\
\quad \quad\quad 
 \sum_{\substack{I\in \bN^n \\ |I|=m }} \frac 1{I!} & H^I(z) 
\frac{\p^m} {\p z^I}
= \sum_{I\in \cC_m}
\frac {(-1)^{\ell (I)} }{ \pi_u( \bar{I}) } \, \, \psi^I \, .
\label{Identity-H(z)-pz-e3}
\end{align}
\end{propo}

\pf First, we assume $o(H(z))\geq 2$ and let 
$F_t(z):=z-tH(z)$ as before. 
By applying $\cSft$ to
Eq.\,(\ref{Psi-Lambda}) with $I=\{m\}$, we have
\begin{align*}
\lambda_m & = (-1)^{m} \sum_{I\in \cC_m } 
\frac{(-1)^{\ell (I)}}{\pi_u(\bar{I})}\,\, \psi^I\, .
\end{align*}
Then, combining the equation above with Eq.\,(\ref{Special-lambda-2}), 
we get Eq.\,(\ref{Identity-H(z)-pz-e1}). 
Furthermore, by Eqs.\,(\ref{Special-psi-m}) and 
(\ref{Special-lambda-3}), 
we get Eqs.\,(\ref{Identity-H(z)-pz-e2}) and 
(\ref{Identity-H(z)-pz-e3}).

Now, consider the case that $o(H(z))<2$. 
Let $u$ be another free variable which is independent with $z$.
Define $\tilde H(z, u):=(u^2H(z), 0)\in 
K\langle\langle z, u \rangle \rangle^{n+1}$.
Then, by applying Eqs.\,(\ref{Identity-H(z)-pz-e1})--(\ref{Identity-H(z)-pz-e3})
 to $\tilde H(z, u)$ 
and then setting $u=1$, we get the identities in the proposition 
for $H(z)$ itself.
\epfv

\subsection{$C_m(z)'s$ in terms of other differential operators in $\Oft$ } 
\label{S5.4}

Motivated by the homogeneous Jacobian conjecture 
that will be discussed in next subsection, below we give  
formulas for $C_m(z)$ $(m\geq 1)$ defined in Lemma \ref{L4.1.1},
which is $(JH)^{m-1} H$
in the commutative case. 

\begin{propo}\label{Cm(z)-others}
For any $m\geq 1$, we have
\begin{align}
C_m(z) &= (-1)^{m} \sum_{I\in \cC_m} (-1)^{\ell(I)} fp (I) \,\, \lambda^I z \, ,  \\
C_m(z) &= (-1)^{m }\sum_{I\in \cC_m } (-1)^{ \ell(I)} lp (I) \,\,  s^I  z  \, ,    \\
C_m(z) &= (-1)^{m }
\sum_{I\in \cC_m} \lp
\sum_{I \succeq J} (-1)^{\ell(J)}  \frac {lp (J)}{ sp (I, J) } \rp \, \phi^I  z  \, , \\
C_m(z) &= \sum_{I\in \cC_m } \lp \sum_{ I \succcurlyeq J  }
(-1)^{\ell(J)- 1} 
\frac {{lp (J)}}{\pi_u(\bar{I}, \bar{J} )} \rp  \, \xi^I  z\, .\label{Cm(z)-others-e4}
\end{align}
In particular, when $z$ are commutative variables, we have
\begin{align*}
(JH)^{m-1} H &= (-1)^{m} \sum_{I\in \cC_m} (-1)^{\ell(I)} fp (I) \,\, \lambda^I  z \, ,  \\
(JH)^{m-1} H  &= (-1)^{m }\sum_{I\in \cC_m } (-1)^{ \ell(I)} lp (I) \,\,  s^I  z  \, ,    \\
(JH)^{m-1} H &= (-1)^{m }
\sum_{I\in \cC_m} \lp
\sum_{I \succeq J} (-1)^{\ell(J)}  \frac {lp (J)}{ sp (I, J) } \rp \, \phi^I  z  \, , \\
(JH)^{m-1} H &= \sum_{I\in \cC_m } \lp \sum_{ I \succcurlyeq J  }
(-1)^{\ell(J)- 1} 
\frac {{lp (J)}}{\pi_u(\bar{I}, \bar{J} )} \rp  \, \xi^I  z\, . 
\end{align*}
\end{propo}

All the formulas above follow directly from the 
identities of $\psi$ with the corresponding differential operators.
For example, by Eq.\,(\ref{Xi-Psi}) with $I=\{m\}$, we have
\begin{align*} 
\Psi_m &=\sum_{K \in \cC_m } \lp \sum_{\bar{K} \succcurlyeq J  }
(-1)^{\ell(J)- 1} 
\frac {{fp (J)}}{\pi_u(\bar{K}, J)} \rp  \, \Xi^K \\
\intertext{Changing the summation index $J$ by $\bar{J}$ and $K$ by $I$:}
 &=\sum_{I\in \cC_m } \lp \sum_{ I \succcurlyeq J  }
(-1)^{\ell(J)- 1} 
\frac {{lp (J)}}{\pi_u(\bar{I}, \bar{J} )} \rp  \, \Xi^I \, .
\end{align*}

Note that, by Eq.\,(\ref{Special-psi-1}),  we have $\psi_m z=C_m(z)$. 
By applying $\cSft$ to the equation above and then 
applying the resulted equation to $z$, 
we get Eq.\,(\ref{Cm(z)-others-e4}).

\subsection{ A Connection of the Jacobian Conjecture with NCSF's}\label{S5.5}

In this subsection, we consider 
the following connection of the well-known 
Jacobian conjecture 
with NCSF's.

Let $K$ be any unital commutative $\bQ$-algebra and
$z=(z_1, z_2, \cdots, z_n)$ be commutative free 
variables. We fix a homogeneous 
$H(z) \in k [z]^{\times n}$ 
of degree $d \geq 2$ and  
$F_t(z)\in \mathbb E^{[d]} 
\langle\langle z \rangle\rangle$, $G_t(z)$ and 
$N_t(z)$ as fixed 
in Section \ref{S3}. 
Denote by $\cD^{[d]}[z]$ the unital 
algebra of the differential 
operators of the polynomial 
algebra $K[z]$, which increase 
the degree by at least $d\geq 2$.
Let $\cS_{F_t}$ be the specialization 
in Theorem \ref{S-Correspondence}.
As one can easily check that, 
for any $F_t(z)\in \mathbb E^{[d]} 
\langle\langle z \rangle\rangle$, 
$\cS_{F_t}$ actually is a 
$K$-Hopf algebra homomorphism from 
$\cNsf$ to $\cD^{[d]}[z]$.
We denote by $\mathcal I_H$ 
the kernel of $\cS_{F_t}$.
Since $\cS_{F_t}$ is a homomorphism 
of $K$-algebras, $\mathcal I_H$ is 
a two-sided ideal of the free 
$K$-algebra $\cNsf$.

By the homogeneous reduction 
in \cite{BCW} and \cite{Y} 
on the Jacobian conjecture, 
it is easy to see that 
the Jacobian conjecture 
is equivalent to 
the following conjecture. 

\vskip2mm

\begin{conj} \label{hgs-JC}
For any homogeneous $H(z)\in K[z]^{\times n}$  
of degree $d \geq 2$, or equivalently, $d=3$, 
let $F_t(z):=z-tH(z)$ and 
$G_t(z):=F_t^{-1}(z)$.
Assume that the Jacobian matrix $JH$ is nilpotent.  
Then the inverse map 
$G_t(z)$ is also a polynomial map 
of $z$ over $K[t]$.
\end{conj}

Note that, by Euler's lemma, 
we have $JH^{m-1}H= \frac 1d (JH)^m  z$. 
Hence the nilpotence of $JH$ implies 
$(JH^{m-1})\,H = 0$ for any $m\geq n$.
By Eqs.\,(\ref{Special-psi-1})
and (\ref{Special-Cm}), this is same 
as saying that the derivations $\psi_m=0$ 
for any $m\geq n$, or equivalently, 
$\Psi_m \in \mathcal I_H$ 
for any $m\geq n$.
On the other hand, by 
Eq.\,(\ref{Special-xi}), 
we see that, $G_t(z)$ is a polynomial map 
iff the derivations $\xi_m=0$ for $m>>0$, 
or equivalently, $\Xi_m \in \mathcal I_H$ 
for $m>>0$. Therefore, by the observations above
and the equivalence of 
the Jacobian conjecture and Conjecture \ref{hgs-JC}, 
we see that the Jacobian conjecture is equivalent to 
the following conjecture.

\begin{conj}\label{JC-NCSF}
For any homogeneous $H(z) \in K[z]^{\times n}$ 
of degree $d=3$, suppose that $\Psi_m \in \mathcal I_H$ 
for any $m\geq n$. 
Then, $\Xi_m 
\in \mathcal I_H$ 
for $m>>0$.
\end{conj}

Note that, by the remarkable symmetric reduction 
on the Jacobian conjecture achieved recently 
in \cite{BE1} and \cite{Me}, 
we may further assume that 
$H(z)$ is the gradient of a homogeneous 
polynomial $P(z)$ of degree $4$, i.e. 
$H(z)=(\frac{\p P}{\p z_1}, \frac{\p P}{\p z_2}, 
\cdots, \frac{\p P}{\p z_n})$.



Therefore, from the point view of 
the Jacobian conjecture, we see
the following open problem becomes 
interesting and important.

\begin{prob}\label{Open-Prob}
For any homogeneous $H(z)\in K[z]^{\times n}$
with the Jacobian matrix 
$JH$ nilpotent, find the relations among NCSF's, 
which decide the two-sided ideal $\mathcal I_H$.
\end{prob}

Note that, by Theorem \ref{StabInjc-best}, 
there are no deciding relations for the ideal 
$\mathcal I_H$, which are independent 
of the choices of $n\geq 1$ and 
$H\in K[z]^{\times n}$. Considering the fact that the classical 
symmetric functions have been well studied and the fact mentioned
in Remark \ref{comm-case}, a good starting point to approach 
Problem \ref{Open-Prob} above might be to consider the case 
when all differential operators in the \cNcs  system $\Oft$ 
commute with each other.

{\small \sc Department of Mathematics, Illinois State University,
Normal, IL 61790-4520.}

{\em E-mail}: wzhao@ilstu.edu.

\end{document}